# CRITICAL THRESHOLDS IN EULER-POISSON EQUATIONS

SHLOMO ENGELBERG, HAILIANG LIU, AND EITAN TADMOR

DEDICATED WITH APPRECIATION TO CIPRIAN FOIAS AND ROGER TEMAM.

ABSTRACT. We present a preliminary study of a new phenomena associated with the Euler-Poisson equations — the so called critical threshold phenomena, where the answer to questions of global smoothness vs. finite time breakdown depends on whether the initial configuration crosses an intrinsic, $O(1)$ critical threshold.

We investigate a class of Euler-Poisson equations, ranging from one-dimensional problem with or without various forcing mechanisms to multi-dimensional isotropic models with geometrical symmetry. These models are shown to admit a critical threshold which is reminiscent of the conditional breakdown of waves on the beach; only waves above certain initial critical threshold experience finite-time breakdown, but otherwise they propagate smoothly. At the same time, the asymptotic long time behavior of the solutions remains the same, independent of crossing these initial thresholds.

A case in point is the simple one-dimensional problem where the unforced inviscid Burgers' solution always forms a shock discontinuity except for the non-generic case of increasing initial profile, $u'_0 \geq 0$. In contrast, we show that the corresponding one dimensional Euler-Poisson equation with zero background has global smooth solutions as long as its initial $(\rho_0, u_0)$- configuration satisfies $u'_0 \geq -\sqrt{2k\rho_0}$ – see (2.11) below, allowing a finite, critical negative velocity gradient. As is typical for such nonlinear convection problems one is led to a Ricatti equation which is balanced here by a forcing acting as a 'nonlinear resonance', and which in turn is responsible for this critical threshold phenomena.

**AMS subject classification:** Primary 35Q35; Secondary 35B30.

**Key Words:** Euler-Poisson equations, critical threshold, breakdown, global existence, isotropic flow.

## Contents

1. INTRODUCTION

The Euler-Poisson equations

$$\rho_t + \nabla \cdot (\rho u) = 0, \quad x \in \mathbb{R}^n, \quad t \in \mathbb{R}^+, \tag{1.1}$$

$$(\rho u)_t + \nabla \cdot (\rho u \otimes u) = k\rho \nabla \phi + \text{viscosity} + \text{relaxation}, \tag{1.2}$$

$$\Delta \phi = \rho + \text{background}, x \in \mathbb{R}^n, \tag{1.3}$$

are the usual statements of the conservation of mass, Newton's second law, and the Poisson equation defining, say, the electric field in terms of the charge. Here $k$ is a given physical constant, which signifies the property of the underlying forcing, repulsive if $k > 0$ and attractive if $k < 0$. The unknowns are the local density $\rho = \rho(x,t)$, the velocity field $u = u(x,t)$, and the potential $\phi = \phi(x,t)$.

This hyperbolic-elliptic coupled system (1.1)-(1.3) describes the dynamic behavior of many important physical flows including charge transport [20], plasma with collision [12], cosmological waves [2] and the expansion of the cold ions [11]. Systems (1.1)- (1.3) also describe the evolution of a star regarded as an ideal gas with self-gravitation ($k < 0$) [16]. The case of repulsive forces ($k > 0$) [22] is relevant for plasma physics. These equations may be obtained from the Vlasov-Poisson- Boltzmann model by setting the mean free path to zero [3]. We would like to point out that the Euler-Poisson equation is closely related to the Schrödinger-Poisson equation via the semi- classical limit and the Vlasov-Poisson equation as well as the Wigner equation. Such relation has been the subject of a considerable number of papers in recent years; we refer to [9], [6] and references therein for further details.

There is a considerable amount of literature available on the global behavior of Euler-Poisson and related problems, from local existence in the small $H^s$-neighborhood of a steady state [16, 21, 8] to global existence of real solutions with geometrical symmetry [5], for the two- carrier types in one dimension [26], the relaxation limit for the weak entropy solution, consult [19] for isentropic case, and [14] for isothermal case.

For the question of global behavior of strong solutions, however, the choice of the initial data and/or damping forces is decisive. The non-existence results in the case of attractive forces have been obtained by Makino-Perthame [18], and for repulsive forces by Perthame [22]. For the study on the singularity formation in the model with diffusion and relaxation, see [27]. In all these cases, the finite life span is due to a *global* condition of large enough initial (generalized) energy, staying outside a critical threshold ball. Using characteristic-based methods, Engelberg [7] gave local conditions for the finite-time loss of smoothness of solutions in Euler-Poisson equations. Global existence due to damping relaxation and with non-zero background can be found in [24, 25, 15]. For the model without damping relaxation the global existence was obtained by Guo [10] assuming the flow is irrotational. His result applies to an $H_2$-small neighborhood of a constant state. Finally we mention the steady solution of non-isentropic Euler-Poisson model analyzed for a collisionless plasma in [17] and for the hydro-dynamic semiconductor in [1] — their approaches are based on the phase plane analysis.

In this paper we present a preliminary study on a new phenomena associated with the Euler-Poisson equations — the so called critical threshold phenomena, where the answer to the question of global vs local existence depends on whether the initial configuration crosses an intrinsic, $O(1)$ critical threshold. Little or no attention has been paid to this remarkable phenomena, and our goal is to bridge the gap between previous studies on the behavior of solutions of the Euler-Poisson equations in the small and in the large. To this end we focus our attention on the n- dimensional isotropic model,

(1.4)
$$r^\nu n_t + (nur^\nu)_r = 0, \quad r > 0,$$
$$\rho(u_t + uu_r) = k\rho\phi_r + \text{viscosity} + \text{relaxation},$$
$$(r^\nu \phi_r)_r = nr^\nu + \text{background}, \quad \nu = n - 1.$$

It is well known that finite time breakdown is a generic phenomena for nonlinear hyperbolic convection equations, which is realized by the formation of shock discontinuities. In the context of Euler-Poisson equation, however, there is a delicate balance between the forcing mechanism (governed by Poisson equation), and the nonlinear focusing (governed by Newton's second law), which supports a critical threshold phenomena. In this paper we show how the persistence of the global features of the solutions hinges on a delicate balance between the nonlinear convection and the forcing mechanism dictated by the Poisson equation as well as other additional forcing mechanism on the right of (1.4).

In particular, the persistence of the global features of solutions does not fall into any particular category (global smooth solution, finite time breakdown, etc.), but instead, these features depend on crossing a critical threshold associated with the initial configuration of underlying problems, very much like the conditional breakdown of waves on the beach; only waves above a certain critical threshold form crests and break down, otherwise they propagate smoothly. See for example [23] for ion-acoustic waves with such critical threshold phenomena.

At present, no rigorous results exist on the question just raised concerning the critical threshold phenomena in equations of Euler-Poisson type. In this paper we provide a detailed account of critical threshold phenomena for a class of Euler-Poisson equations without pressure forcing. We use these systems to demonstrate the ubiquity of *critical thresholds* in the solutions of some of the equations of mathematical physics. At the same time, we show that the asymptotic long time behavior of the solutions remains the same, independent of whether the initial data has crossed the critical threshold or not.

We note in passing that in this paper we restrict ourselves to the pressureless isotropic Euler-Poisson equations. The existence of the pressure allows for additional balance, and we hope to explore the critical threshold phenomena for the general, possibly non-isotropic model with additional pressure forces in the future work.

A simple example of an equations with a critical threshold is the 1D unforced Burgers' equation, $u_t + uu_x = 0$. This equation describes the movement of particles that are not being acted on by any forces. The variable $u(x,t)$ represents the velocity of the particle located at position $x$ at time $t$. The global existence is ensured if and only if $u_0' \geq 0$. Thus, the Burgers' solution forms a shock discontinuity unless its initial profile is monotonically non-decreasing. In this case the finite time breakdown of the Burgers' solution is a generic phenomena. In contrast we show below that the corresponding 1D Euler-Poisson equation with zero background has global smooth solutions if and only if $u_0' > -\sqrt{2k\rho_0}$ (see (2.11) below), allowing a finite negative velocity gradient. This is the critical threshold phenomena we are referring to. As is typical for such nonlinear convection problems, one is led to a general Ricatti equation which is complemented by a particular inhomogeneous forcing dictated by the Poisson equation. It is the delicate balance of the latter, acting as a nonlinear resonance, which is responsible for this critical threshold phenomena.

This paper is organized as follows. In Section 2, we consider the 1D Euler-Poisson equations with zero background, and show that the solutions of the corresponding Cauchy problem blow up in finite time if and only if certain local "threshold" conditions on the initial data are met. In this case the density and the velocity gradient are shown to decay at some algebraic rates. For this simple model we utilize both the Eulerian and the Lagrangian description of the flows to investigate the critical threshold phenomena. We also discuss how the solution behavior depends on the initial data as well as on the coupling parameter $k$. The behavior of the velocity gradient is outlined for various cases classified by the relative size of the initial data.

In Section 3, we present the critical thresholds for the 1D model with nonzero constant background and the effects of the damping relaxation. In this case the solution oscillates. The oscillatory behavior of the solution is induced by the presence of nonzero background, and the oscillations can be prevented only by strong damping relaxation. The zero limit of the background is shown to be a kind of "singular limit", which coincides with (close to) the zero background case studied in Section 2. We also justify a rather remarkable phenomena, namely, that the non-zero background is able to balance both the nonlinear convective focusing effects and the attractive forces, to form a global smooth solution subject to a critical threshold in the attractive case $k < 0$. Moreover in this case the density converges to the vacuum exponentially fast in time.

In Section 4 we continue the discussion for the 1D model with viscosity, which is similar (though not identical) to the one found in the Navier-Stokes equations. We obtain an upper threshold for the existence of global smooth solution and an lower threshold for the finite time breakdown, and consequently, these imply the existence of a critical threshold. If the initial data happens to be below the lower threshold, the solution must breakdown in finite time. We see that the presence of a self-induced electric field as well as additional viscosity are not necessarily enough to stop the formation of singularities.

In Section 5 we consider the Euler-Poisson equations governing the $\nu + 1$ dimensional isotropic ion expansion in the electrostatic fluid approximation for cold ions. We show that for each model with given integer $\nu > 0$, there exists a critical threshold condition. Several issues which are also clarified in this section include:

(1) Expansion rate of the flow path (consult the recent work [6] of Dolbeault and Rein in this context).

(2) The large time behavior of the velocity;

(3) The decay rate of the density as well as the velocity gradient;

(4) Sharp estimate of blow up time when the initial data exceed the critical threshold.

More precisely we provide the explicit form of the critical threshold for the planar case and 4-dimensional case. For other cases (including the cylindrical and the spherical case) we confirm the critical threshold phenomena by establishing both the upper threshold for the existence of the global smooth solution and the lower threshold for the finite time breakdown. A key step in the proof is to introduce the proper weighted electric field, which is shown to be constant along the particle path. This fact combined with the momentum equation gives the decoupled equation for the flow map, and a nonlinear resonance is responsible for the critical threshold phenomena.

## 2. Critical thresholds–1D model with zero background

We consider the 1D Euler-Poisson equation of the form

$$\rho_t + (\rho u)_x = 0, \tag{2.1}$$

$$u_t + u u_x = -k\phi_x = kE, \tag{2.2}$$

$$E_x = -\phi_{xx} = \rho. \tag{2.3}$$

Here $k$ is a given physical constant, which signifies that the underlying forcing $kE$ is repulsive if $k > 0$ and attractive if $k < 0$. The unknowns are the density of negatively charged matter $\rho = \rho(x,t)$, the velocity field $u = u(x,t)$, the electric field $E = E(x,t)$, and the potential $\phi = \phi(x,t)$. Equation (2.1) is a statement of the conservation of mass, equation (2.2) is a statement of Newton's second law, and equation (2.3) defines the electric field in terms of the charge.

To solve the problem on the half plane $(x,t) \in \mathbb{R} \times \mathbb{R}^+$ we prescribe initial data as

(2.4) $$\rho(x,0) = \rho_0(x) > 0, \quad \rho_0 \in C^1(\mathbb{R}),$$

(2.5) $$u(x,0) = u_0(x), \quad u_0 \in C^1(\mathbb{R}),$$

and we show that the solutions of (2.1)-(2.3) with the above initial data break down in finite-time if and only if certain local "threshold" conditions on the initial- data are met.

Set $d := u_x(x,t)$, then $\partial_x$(2.2) together with (2.1) yield by differentiation along the characteristics,

$$d' + d^2 = k\rho,$$
$$\rho' + \rho d = 0, \quad ' := \partial_t + u\partial_x.$$

Multiply the first equation by $\rho$, the second equation by $d$ and take the difference. This gives

$$\left(\frac{d}{\rho}\right)' = \frac{\rho d' - d\rho'}{\rho^2} = k,$$

and upon integration one gets

$$\frac{d}{\rho} = \beta(t) \quad \text{with} \quad \beta(t) := kt + u_0'/\rho_0.$$

The decoupled equations for $d$ and $\rho$ now read

(2.6) $$d' + d^2 = \frac{k}{\beta(t)} d,$$

(2.7) $$\rho' + \beta(t)\rho^2 = 0.$$

¿From these equations we can obtain the explicit solution formula for $d = u_x$ and $\rho$, respectively. We want to point out that it is this time-dependent factor $\beta(t)$ balancing the nonlinear quadratic term that is responsible for the critical threshold phenomena. Indeed, with $\frac{k}{\beta(t)} \equiv 0$ one has the usual blow up, $d(t) = \frac{d_0}{1+d_0 t}$, associated with the unforced Ricatti equation $d' + d^2 = 0$. The presence of a forcing of similar strength on the right of (2.6), $\frac{k}{kt+u_0'/\rho_0}d$, leads to 'nonlinear resonance', which, as we shall see below, prevents blow up at least above a critical threshold, such that $-\sqrt{2k\rho_0} < u_0' < 0$. To highlight this fact we present the following general lemma.

**Lemma 2.1.** *Consider the ODE*

$$w_t = a(t)w + b(t)w^2, \quad w(t=0) = w_0.$$

*It admits a global solution*

$$w(t) = \frac{w_0 e^{\int_0^t a(\tau)d\tau}}{1 - w_0 \int_0^t B(\tau)d\tau}, \quad B(t) := b(t)e^{\int_0^t a(\tau)d\tau}.$$

*provided the initial data, $w_0$, is prescribed so that*

$$w_0 \int_0^t B(\tau)d\tau < 1 \quad \text{for all} \quad t > 0.$$

*Proof.* Set
$$v(t) = w(t)e^{-\int_0^t a(\tau)d\tau}.$$
Substitution into the above ODE leads to
$$v_t = B(t)v^2.$$
Note that $v_0 = w_0$. Then the solution can be written explicitly as
$$v = w_0\Big(1 - w_0 \int_0^t B(\tau)d\tau\Big)^{-1},$$
from which the lemma immediately follows. □

As an immediate application of this lemma we check the conditions for $d_0$ and $\rho_0$ so that the global existence of $u_x$ and $\rho$ is ensured. Lemma 2.1 applies to the above equation for $d$, (2.6), with $b(t) \equiv -1$, $a(t) = \frac{k}{\beta(t)}$ and $w_0 = d_0 = u'_0$. It follows that if $u'_0 > -\sqrt{2k\rho_0}$ then global regularity for $d$ is ensured, for
$$1 - w_0 \int_0^t B(\tau)d\tau = 1 + u'_0 \int_0^t e^{\int_0^\tau \frac{k}{\beta(s)}ds} d\tau d\tau = 1 + u'_0 t + \frac{k}{2}\rho_0 t^2 > 0.$$
Similarly, for the $\rho$ equation (2.7) we have $a \equiv 0$, $b = -\beta$ and $w_0 = \rho_0$, and hence
$$1 - w_0 \int_0^t B(\tau)d\tau = 1 + \rho_0 \int_0^t \beta(\tau)d\tau = 1 + u'_0 t + \frac{k}{2}\rho_0 t^2 > 0,$$
provided the initial velocity gradient $u'_0$ remains above the same critical threshold $-\sqrt{2k\rho_0}$.

We now turn to an alternative derivation of this critical threshold. This Lagrangian-like approach will prove to be useful for more general cases. We start by appealing to "physical" considerations. As the system being described has no external forces acting on it, momentum ought to be conserved. In fact, multiplying (2.2) by $\rho$, multiplying (2.1) by $u$, and adding the resulting equations, one finds the momentum equation in its standard conservative form
$$(\rho u)_t + (\rho u^2)_x = k\rho E = (kE^2/2)_x, \quad \rho = E_x.$$
It follows that the momentum is conserved, $\int_{-\infty}^{\infty} \rho u(\cdot, t)dx =$ const, provided the boundary terms vanish, i.e., $\rho u^2 \to 0$ as $x \to \pm\infty$, and in particular, $E^2(\infty) - E^2(-\infty) = 0$. It is reasonable to require that the total charge, $E(\infty)$, be finite, and since $E_x = \rho \geq 0$, this implies that $E(\infty) = -E(-\infty)$, for otherwise, $\rho \equiv 0$. Thus the electric field is given by
$$E(x,t) = \frac{1}{2}\left(\int_{-\infty}^x \rho(\xi, t)\,d\xi - \int_x^\infty \rho(\xi, t)\,d\xi\right).$$

Equipped with this expression of $E$ in terms of the density $\rho$, we can obtain the explicit solution along the characteristic curves, $x(\alpha, t)$, parameterized with respect to the initial positions, $x(\alpha, 0) = \alpha$,

(2.8) $$\frac{d}{dt}x(\alpha, t) = u(x(\alpha, t)), \quad x(\alpha, 0) = \alpha.$$

The momentum equation (2.2) tells us that

(2.9) $$\frac{d}{dt}u(x(\alpha, t), t) = kE(x(\alpha, t), t),$$

and hence the electric field remains constant along $x(\alpha, t)$,

$$\begin{aligned} 2\frac{d}{dt}E(x(\alpha,t),t) &= 2\frac{d}{dt}(x(\alpha,t)) \cdot \rho(x(\alpha,t),t) + \int_{-\infty}^{x(\alpha,t)} \rho_t(\xi,t)\,d\xi - \int_{x(\alpha,t)}^{\infty} \rho_t(\xi,t)\,d\xi \\ &= 2u(x(\alpha,t),t)\rho(x(\alpha,t),t) \\ &\quad - \int_{-\infty}^{x(\alpha,t)} (\rho(\xi,t)u(\xi,t))_x\,d\xi + \int_{x(\alpha,t)}^{\infty} (\rho(\xi,t)u(\xi,t))_x\,d\xi \\ &= 0. \end{aligned}$$

Physically, the constancy of the electric field along characteristics is clear: as no charge can cross trajectories, the amount of charge to the right and to the left of any given trajectory is constant. And since the electric field along a trajectory is half of the difference of these numbers, the electric field on any trajectory must be constant as well. With $E(\alpha, 0) =: E_0, u(\alpha, 0) =: u_0$, and $\rho(\alpha, 0) =: \rho_0$, we find from (2.9) that

$$u(x(\alpha,t)) = u_0 + kE_0 t.$$

This together with (2.8) yield

(2.10) $$x(\alpha, t) = \alpha + u_0 t + kE_0 t^2/2.$$

In the inviscid Burgers' equation (corresponding to the case $k = 0$) the straight characteristics must intersect in finite time leading to finite time breakdown. Here the straight characteristic curves are replaced by the characteristic parabolas, which explain the critical threshold phenomena. Indeed, since $u = u_0 + kE_0 t$ and $E_{0\alpha} = \rho_0$, we conclude

$$u_x(x(\alpha,t),t) = (u_0' + k\rho_0 t)/\frac{\partial x(\alpha,t)}{\partial \alpha} = \frac{u_0' + k\rho_0 t}{1 + u_0' t + k\rho_0 t^2/2}, \quad u_0' := \frac{\partial u_0(\alpha)}{\partial \alpha}.$$

Integrating the $\rho$-equation, (2.1), which we rewrite as $\frac{d}{dt}\rho(x(\alpha,t),t) = -u_x\rho$, we find that

$$\rho(x(\alpha,t),t) = \frac{\rho_0}{\Gamma(\alpha,t)}, \quad \Gamma(\alpha,t) := 1 + u_0't + k\rho_0 t^2/2.$$

Clearly, the positivity of the so called "indicator" function, $\Gamma(\alpha, t) := \frac{\partial x(\alpha,t)}{\partial \alpha}$, is a necessary and sufficient condition for existence of global smooth solution, $|u_x|, \rho \leq \text{Const.}$, and consequently, higher derivatives are bounded. Thus the solution remains smooth for all time if and only if $-u_0' < \sqrt{2k\rho_0}$. Conversely, if there are points at which this condition is not fulfilled, then $\rho$ and $u_x$ will blow up at finite-time.

To sum up, we state

**Theorem 2.2.** *The system of Euler-Poisson equations (2.1)-(2.3) admits a global smooth solution if and only if*

(2.11) $$u_0'(\alpha) > -\sqrt{2k\rho_0(\alpha)}, \quad \forall \alpha \in \mathbb{R}.$$

*In this case the solution of (2.1)- (2.3) is given by*

$$\rho(x(\alpha,t),t) = \frac{\rho_0}{\Gamma(\alpha,t)}, \quad u_x(x(\alpha,t),t) = \frac{u_0' + k\rho_0 t}{\Gamma(\alpha,t)}, \quad \Gamma(\alpha,t) := 1 + u_0't + k\rho_0 t^2/2$$

*so that $\rho \sim t^{-2}$ and $u_x \sim t^{-1}$ as long as $\rho_0 \neq 0$. If condition (2.11) fails, then the solution breaks down at the finite time, $t_c$, where $\Gamma(\alpha, t_c) = 0$.*

To gain further insight into the behavior of the solution, we now turn to discuss the dependence of the solution on the relative size of $d_0 = u_0'$ and $\rho_0$ as well as the parameter $k$.

To be specific we consider only the behavior of $d \equiv u_x$ for the repulsive forces $k > 0$, since the solution always breaks down in the attractive case $k < 0$. We first look at the dependence on $d_0 = u'_0$. Figures 1-3 describe the three different scenarios for the evolution of $u_x$ depending on the relative size of $d_0$ and $\rho_0$.

(1) $d_0 > 0$ (see Figure 1). There are two such cases:
    (i) if $d_0 > \sqrt{k\rho_0}$, then $d$ is decreasing and satisfies

$$0 < d \leq d_0, \quad d \sim 2/t.$$

    (ii) if $0 < d_0 < \sqrt{k\rho_0}$, then we have

$$0 < d \leq d_{\max}, \quad d_{\max} = \frac{k\rho_0}{\sqrt{2k\rho_0 - d_0^2}},$$

where $d_{\max}$ denotes its local maximum taken at $t_e^+$. At this time $d' = 0$, therefore $t_e$ satisfies

$$\frac{k^2 \rho_0^2}{2} t_e^2 + k\rho_0 d_0 + d_0^2 - k\rho_0 = 0, \quad d_{max} = d(t_e^+).$$

(2) $-\sqrt{2k\rho_0} < d_0 < 0$ (see Figure 2). There are also two cases:
    (i) if $-\sqrt{k\rho_0} < d_0 < 0$, then $d$ starts to increase and becomes zero at $t_0 = -\frac{d_0}{k\rho_0}$ and then attains its maximum at

$$t_e^+ = \frac{\sqrt{2k\rho_0 - d_0^2} - d_0}{k\rho_0} > t_0.$$

In this case we have

$$d_0 \leq d \leq d_{\max} = \frac{k\rho_0}{\sqrt{2k\rho_0 - d_0^2}}, \quad d \sim 2/t.$$

    (ii) $-\sqrt{2k\rho_0} < d_0 < -\sqrt{k\rho_0}$. In this case we have

$$d_{\min} \leq d \leq d_{max}, \quad d \sim 2/t,$$

where

$$d_{\min} = d(x, t_e^-) = \frac{-k\rho_0}{\sqrt{2k\rho_0 - d_0^2}} \quad d_{\max} = d(x, t_e^+) = \frac{k\rho_0}{\sqrt{2k\rho_0 - d_0^2}}$$

and

$$t_e^\pm = \frac{\pm\sqrt{2k\rho_0 - d_0^2} - d_0}{k\rho_0}.$$

(3) $d_0 < -\sqrt{2k\rho_0}$ (see Figure 3). In this case the solution must break down at time $t = t_c^-$. The blow up time $t_c$ can be obtained via $\Gamma(\alpha, t_c^\pm) = 0$, that is,

$$t_c^\pm = \frac{-d_0 \pm \sqrt{d_0^2 - 2k\rho_0}}{k\rho_0}.$$

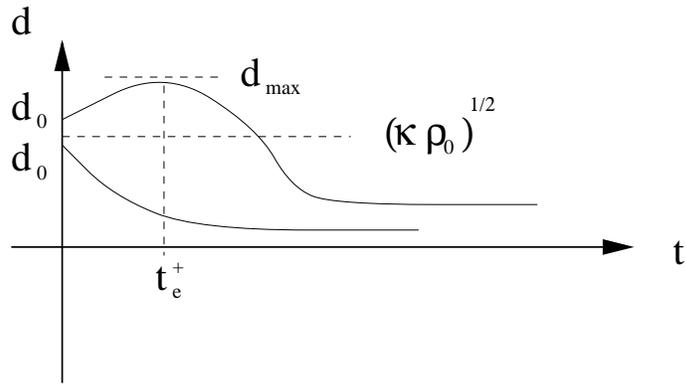

FIGURE 1. $0 < d_0$, $k > 0$

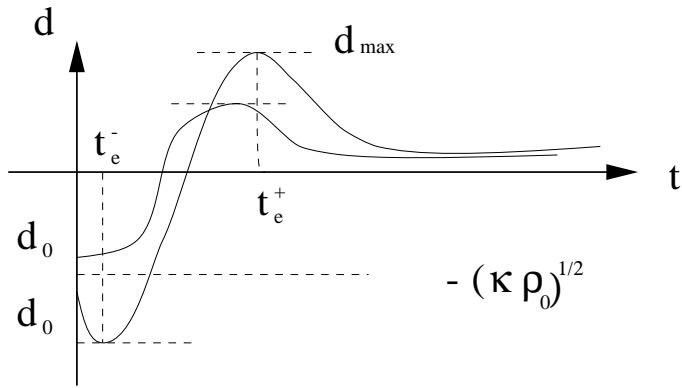

FIGURE 2. $-\sqrt{2k\rho_0} < d_0 < 0$, $k > 0$

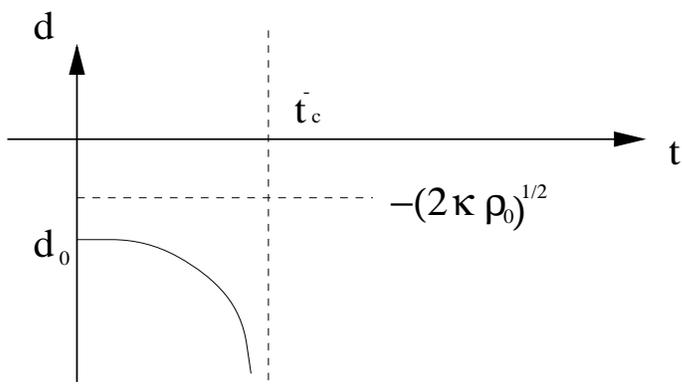

FIGURE 3. $d_0 < -\sqrt{2k\rho_0}$, $k > 0$

Next we look at the solution behavior as the parameter $k$ changes. Figures 4-7 display such changes for various choices of $k$. Again there are three possible scenarios:

(1) $k > d_0^2/\rho_0$ (see Figure 4). There are two cases:
   (i) if $d_0 > 0$, then $0 < d \leq \frac{k\rho_0}{\sqrt{2k\rho_0 - d_0^2}}$;
   (ii) if $d_0 < 0$, then $d_0 \leq \frac{k\rho_0}{\sqrt{2k\rho_0 - d_0^2}}$.

(2) $\frac{d_0^2}{2\rho_0} < k < d_0^2/\rho_0$ (see Figure 5)
   (i) if $d_0 > 0$, then $0 < d \leq d_0$;
   (ii) if $d_0 < 0$, then $\frac{-k\rho_0}{\sqrt{2k\rho_0 - d_0^2}} \leq d \leq \frac{k\rho_0}{\sqrt{2k\rho_0 - d_0^2}}$.

(3) $0 < k < \frac{d_0^2}{2\rho_0}$ (see Figure 6).
   (i) if $d_0 > 0$, then $d$ is decreasing and $0 < d \leq d_0$.
   (ii) if $d_0 < 0$, then $d \leq d_0$ and $d$ starts to decrease and becomes unbounded at time $t_c^-$.

In closing, let us note the remaining cases of $k$: if $k = 0$ then we have the decoupled Burgers' equation; and if $k < 0$, the solution always breaks down. In either case, there is no critical threshold which yields global smooth solution for an admissible set of initial data.

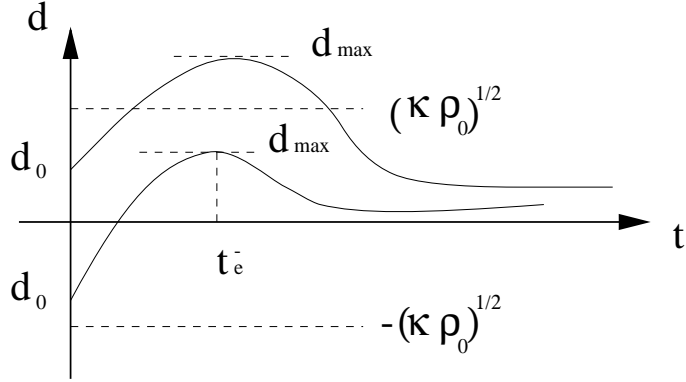

FIGURE 4. $k > d_0^2/\rho_0$

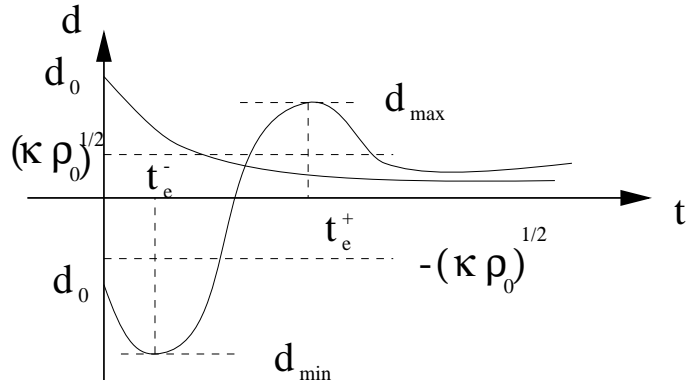

FIGURE 5. $\frac{d_0^2}{2\rho_0} < k < \frac{d_0^2}{\rho_0}$

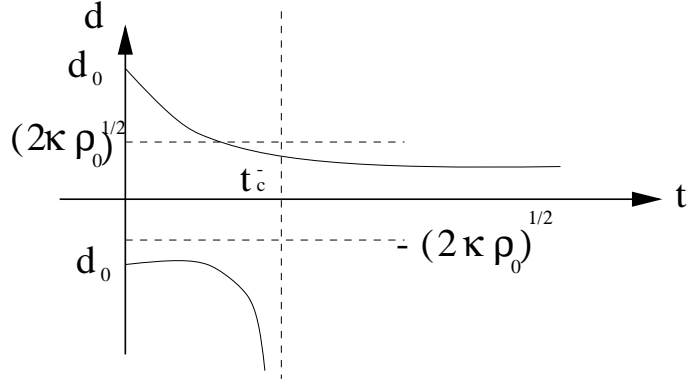

FIGURE 6. $0 \leq k < \frac{d_0^2}{2\rho_0}$

## 3. CRITICAL THRESHOLDS– 1D MODELS WITH NONZERO BACKGROUND

### 3.1. The basic model with constant background.
We now consider the system

(3.1) $$\rho_t + (\rho u)_x = 0,$$

(3.2) $$u_t + uu_x = -k\phi_x = kE,$$

(3.3) $$E_x = -\phi_{xx} = \rho - c,$$

with constant "background" state $c > 0$. Here we require that $\int_{-\infty}^{\infty}(\rho(\xi) - c)\,d\xi = 0$. The presence of a constant background $c$ in the Poisson equation changes the "physical situation:" we are now working in a universe which has a fixed background charge density of $-c$. There is also an equal amount of movable charge, $\rho(x)$. This is (an approximation of) the situation, for example, inside a metal or a doped semiconductor. The fixed background charge of $-c$ corresponds to the fixed positive charge of an element which has had an electron stripped from its outermost shell. The movable charge corresponds to the electrons that have been liberated from the atoms.

Using the Lagrangian-like approach, we solve the system of equations (3.1)-(3.3). For this system, one cannot expect the total momentum of the negatively charged particles to be conserved–they are being acted on by an outside force. From (3.3) we find that $E(x) = \int_{-\infty}^{x}(\rho(\xi) - c)\,d\xi$. As the net charge in our universe is zero, we expect that the electric field intensity to vanish at $x = \pm\infty$, and hence we require that $\rho(\pm\infty, t) = c$. Likewise we require the particles be at rest at far field, i.e., $u(\pm\infty, t) = 0$. This says that far from the origin our system is "properly charge balanced" and at rest.

As noted above, (3.2) says that $\frac{d}{dt}u = kE$ with $\frac{d}{dt}$ denoting the usual differentiation along the characteristics, $\frac{d}{dt}x(\alpha, t) = u(x(\alpha, t), t)$. Using $E(x) = \int_{-\infty}^{x}(\rho(\xi) - c)\,d\xi$, and following the same basic steps as above, we find that $\frac{d}{dt}E = -cu$. Combining these two results, we arrive at

(3.4) $$\frac{d^2u}{dt^2} = -cku, \quad u(\alpha, 0) = u_0, \quad u_t(\alpha, 0) = kE_0,$$

yielding, for $k > 0$,

$$u(x(\alpha, t), t) = u_0 \cos(\sqrt{ck}t) + \frac{kE_0}{\sqrt{ck}} \sin(\sqrt{ck}t).$$

Here, parabolas (corresponding to $c = 0$) are replaced by a different geometry of characteristics, where $\frac{d}{dt}x(\alpha,t) = u(x(\alpha,t),t)$ implies

$$x(\alpha,t) = \frac{u_0}{\sqrt{ck}}\sin(\sqrt{ck}t) + \frac{E_0}{c}(1 - \cos(\sqrt{ck}t)) + \alpha.$$

Expressed in terms of the indicator function, $\Gamma(\alpha,t) := x_\alpha(\alpha,t)$, given by

$$\Gamma(\alpha,t) = 1 + \frac{u_0'(\alpha)}{\sqrt{ck}}\sin(\sqrt{ck}t) - \frac{\rho_0(\alpha) - c}{c}(\cos(\sqrt{ck}t) - 1),$$

we proceed as before to find that $u_x = u_\alpha/x_\alpha$ is given by

(3.5) $$u_x(x(\alpha,t),t) = \frac{u_0'(\alpha)\cos(\sqrt{ck}t) + k\frac{\rho_0(\alpha)-c}{\sqrt{ck}}\sin(\sqrt{ck}t)}{\Gamma(\alpha,t)}.$$

Note that $u_x = \Gamma_t/\Gamma$ we find from (3.1) that $\rho' = -\rho u_x = -\rho\Gamma_t/\Gamma$ which in turn leads to

(3.6) $$\rho(x(\alpha,t),t) = \frac{\rho_0(\alpha)}{\Gamma(\alpha,t)}.$$

Clearly, there is a global smooth solution if and only if $\Gamma(\alpha,t)$ remains positive. For this to hold, we note that there exists $\tau$ such that $\Gamma$ can be rewritten as,

$$\Gamma(\alpha,t) = \frac{\rho_0(\alpha)}{c} + \sqrt{\frac{u_0'(\alpha)^2}{ck} + \left(\frac{\rho_0(\alpha)}{c} - 1\right)^2}\sin(t + \tau),$$

and hence, $\Gamma(\alpha,t) > 0$ if and only if

$$\sqrt{\frac{u_0'(\alpha)^2}{ck} + \left(\frac{\rho_0(\alpha)}{c} - 1\right)^2} < \frac{\rho_0(\alpha)}{c}.$$

This is equivalent to the condition that $|u_0'(\alpha)| < \sqrt{k(2\rho_0(\alpha) - c)}$ for all $\alpha \in \mathbb{R}$.

We can summarize the case with repulsive force $k > 0$

**Theorem 3.1.** *Consider the system of Euler-Poisson equations (3.1)-(3.3) with constant background charge c and the repulsive force $k > 0$. Then it admits a global smooth solution if and only if*

(3.7) $$|u_0'(\alpha)| < \sqrt{k(2\rho_0(\alpha) - c)}, \quad \forall \alpha \in \mathbb{R}.$$

*In this case, the density oscillates around the nonzero background charge c, and the velocity gradient does not decay in time. If condition (3.7) fails, however, the solution breaks down at finite time.*

*Remark* 3.1. The above threshold is sharp, as shown in the phase plane $(\Gamma, \Gamma_t)$. Actually integration of (3.4) along characteristic curves $x(\alpha,t)$ yields

(3.8) $$x'' + ckx = kE_0(\alpha) + ck\alpha,$$

and differentiation with respect to $\alpha$, combined with $E_0' = \rho_0 - c$, leads to

$$\Gamma'' + ck\Gamma = k\rho_0.$$

Its energy integral then becomes

$$\frac{1}{2}\Gamma'^2 + ck\Gamma^2 - k\rho_0\Gamma = \frac{1}{2}u_0'^2 - k(\rho_0 - c),$$

with a trajectory which is an ellipse centered at $(\rho_0/c, 0)$. Here, even if initially $u'_0 > 0$, then after some time $\Gamma_t$ can still become negative, with the possible breakdown due to intersection with the $\Gamma = 0$ line. Condition (3.7) is the precise condition which rules out such a scenario. Geometrically, as $c$ tends to zero, the center of the above ellipse moves to the infinity and the closed curve splits into the parabola we met earlier with the zero background case. In this limiting case, once $u'_0 > 0$, the trajectory will always run away from the 'singular' axis $\Gamma = 0$, see Figures 7-8.

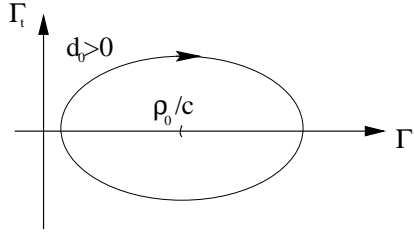

FIGURE 7. $c > 0$

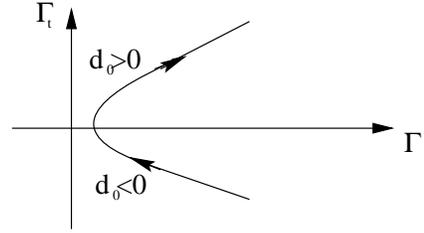

FIGURE 8. $c = 0$

*Remark 3.2.* Note that as $c \to 0$, the critical threshold condition (3.7) tends to $|u'_0(\alpha)| < \sqrt{2k\rho_0(\alpha)}$. This is almost the same condition we found when we considered (2.1)-(2.3), the condition given in (2.11). Here, however, the part of the condition that involves $u'_0(\alpha)$ appears with an absolute value even in the limit as $c \to 0$. The absolute value appears because of the oscillatory nature of the solutions of our equation. When the derivative of the initial condition is too negative for some large value of $t$, the equation loses smoothness. One finds that as $c \to 0$, this time tends to infinity. The passage to the case $c = 0$ is, therefore, a kind of "singular limit."

We now turn to the case of an attractive force, $k < 0$. Using the same basic steps as before we get the same characteristic description, (3.8) as the case $k > 0$. Integration of this equation subject to initial data $(x, x')(\alpha, 0) = (\alpha, u_0(\alpha))$ yields

$$x(\alpha, t) = \alpha + \frac{E_0}{c} + \frac{1}{2}\left(-\frac{E_0}{c} - \frac{u_0}{\sqrt{-ck}}\right) e^{-\sqrt{-ck}t} + \frac{1}{2}\left(-\frac{E_0}{c} + \frac{u_0}{\sqrt{-ck}}\right) e^{\sqrt{-ck}t}.$$

Therefore, using $E_{0\alpha} = \rho_0(\alpha) - c$, we find the indicator function $\Gamma(\alpha, t) := \frac{\partial x(\alpha,t)}{\partial \alpha}$ given by

$$(3.9) \quad \Gamma(\alpha, t) = \frac{\rho_0}{c} + \frac{1}{2}\left(1 - \frac{\rho_0}{c} - \frac{u'_0}{\sqrt{-ck}}\right) e^{-\sqrt{-ck}t} + \frac{1}{2}\left(1 - \frac{\rho_0}{c} + \frac{u'_0}{\sqrt{-ck}}\right) e^{\sqrt{-ck}t}.$$

We conclude the rather remarkable phenomena, namely that the non-zero background is able to balance both the nonlinear convective focusing effects as well as the attractive forces, to form a global smooth solution subject to a critical threshold.

**Theorem 3.2.** *Consider the system of Euler-Poisson equations (3.1)-(3.3) with constant background charge $c > 0$ and subject to attractive force, $k < 0$. Then, it admits a global smooth solution if and only if*

$$(3.10) \quad u'_0(\alpha) \geq -\left(1 - \frac{\rho_0(\alpha)}{c}\right)\sqrt{-ck}, \quad \forall \alpha \in \mathbb{R}.$$

*In this case, the density approaches the zero exponentially in time, and the velocity gradient remains bounded uniformly in time. If condition (3.10) fails, then the solution breaks down in finite time.*

*Proof.* As argued in the proof of Theorem 3.1, we have that

$$\rho(x(\alpha,t),t) = \rho_0(\alpha)/\Gamma(\alpha,t), \quad u_x(x(\alpha,t),t) = \Gamma_t(\alpha,t)/\Gamma(\alpha,t).$$

It is necessary and sufficient to show that the indicator function $\Gamma(\alpha,t) > 0$ for all $t > 0$ if and only if (3.10) is met. The necessity is obvious since otherwise if (3.10) fails then the second parenthesis on the right of (3.9) is negative, and hence $\Gamma(\alpha,t)$ would become negative for $t$ large.

For the sufficiency, there are two cases: either $u_0'(\alpha) = -(1 - \frac{\rho_0(\alpha)}{c})\sqrt{-ck}$, in which case

$$\Gamma(\alpha,t) = \frac{\rho_0}{c} + \left(1 - \frac{\rho_0}{c}\right)e^{-\sqrt{-ck}t},$$

and obviously, $\Gamma(\alpha,t)$ remains positive for all time $t > 0$. For the other case where $u_0' > -\left(1 - \frac{\rho_0(\alpha)}{c}\right)\sqrt{-ck}$, the possible zeros $t^*$ of $\Gamma(\alpha,t)$ are determined by

$$(3.11) \qquad e^{\sqrt{-ck}t^*} = \frac{-\frac{\rho_0}{c} + \sqrt{\frac{2\rho_0}{c} - 1 + \frac{u_0'^2}{-ck}}}{1 - \frac{\rho_0}{c} + \frac{u_0'}{\sqrt{-ck}}}.$$

A simple check shows that the quantity on the right is less than 1, that is, either $t^* < 0$ or no such real $t^*$ exists. Therefore $\Gamma(\alpha,t) > 0$ for all time $t > 0$. $\square$

*Remark* 3.3. Let us, again, note that with the zero background model, the repulsive force, $k > 0$, is necessary for the global existence of the smooth solution. When the nonzero background is being taken into account, we could still have the global existence even when the force is attractive, $k < 0$. The balancing effect of $k$ and $c$ can be observed clearly from the above results. However, in both cases, we find that there is a local " critical threshold" condition on the initial data such that the solution remains smooth for all time if and only if this condition is met.

**3.2. A constant background model with relaxation.** We consider a further modification of our problem (3.1), (3.3), where (3.2) is now augmented by a relaxation term

$$(3.12) \qquad u_t + uu_x = -k\phi_x - \frac{u}{\epsilon} = kE - \frac{u}{\epsilon}, \quad \epsilon > 0.$$

We still require that $\int_{-\infty}^{\infty}(\rho(\xi) - c)\,d\xi = 0$.

We are now working in a universe which has a fixed background charge density of $-c$. There is also an equal amount of movable charge, $\rho(x)$. This is (an approximation of) the situation inside a metal or a doped semiconductor. The term $-u/\epsilon$ is a "friction term", which, as we shall see, causes solutions to decay.

As before, we cannot expect the total momentum of the negatively charged particles to be conserved; after all, they are being acted on by an outside force. As the net charge in our universe is zero, we expect that the electric field intensity at $x = \pm\infty$ will be zero. Consequently we require that $\rho(\pm\infty, t) = c$ and $u(\pm\infty, t) = 0$ so that $E(x) = \int_{-\infty}^{x}(\rho(\xi) - c)\,d\xi \to 0$ as $x \to \infty$.

Using the same techniques that we have used previously, we find that as a consequence of (3.12) we have,

$$\frac{d}{dt}u = kE - u/\epsilon,$$

and proceeding as before,

$$\frac{d^2u}{dt^2} = -cku - \frac{1}{\epsilon}\frac{d}{dt}u.$$

Let $\epsilon > 1/\sqrt{4ck}$, and define $\mu = \sqrt{ck - \frac{1}{4\epsilon^2}}$. This guarantees that the solution of the last ODE will consist of damped sinusoids and decaying exponentials, where we find,

$$u(x(\alpha,t),t) = e^{\frac{-t}{2\epsilon}}\left(u_0\cos(\mu t) + \frac{kE_0 - u_0/(2\epsilon)}{\mu}\sin(\mu t)\right).$$

We see that as long as the solution of our system is well defined, $u(x,t)$ must vanish as $t \to \infty$, and since $\frac{d}{dt}x(\alpha,t) = u(x(\alpha,t))$, we conclude that the characteristic curves have the form

$$x(\alpha,t) = \alpha + E_0/c + \frac{e^{\frac{-t}{2\epsilon}}}{c}\left(-E_0\cos(\mu t) + \frac{2\epsilon cu_0 - E_0}{2\epsilon\mu}\sin(\mu t)\right).$$

As we have seen already, solutions will cease to exist when the indicator function $\Gamma(\alpha,t) = \frac{\partial}{\partial\alpha}x(\alpha,t)$ vanishes, for

$$u_x(x(\alpha,t),t) = \frac{\partial}{\partial\alpha}u(x(\alpha,t),t)\frac{1}{\Gamma(\alpha,t)} \uparrow \infty.$$

It is easy to see that this leads to a critical threshold condition, namely that the vanishing of the following indicator function

$$(3.13) \qquad \Gamma(\alpha,t) = \frac{\rho_0}{c} + \frac{e^{\frac{-t}{2\epsilon}}}{c}\left((c-\rho_0)\cos(\mu t) + \frac{2\epsilon cu'_0 - \rho_0 + c}{2\epsilon\mu}\sin(\mu t)\right).$$

To make precise the condition on the initial data to have the global solution, we need to verify the first time when the local minimum of $\Gamma(\alpha,t) = \frac{\partial}{\partial\alpha}x(\alpha,t)$ intersect the 'singular' $\Gamma = 0$ axis as shown in Figure 10.

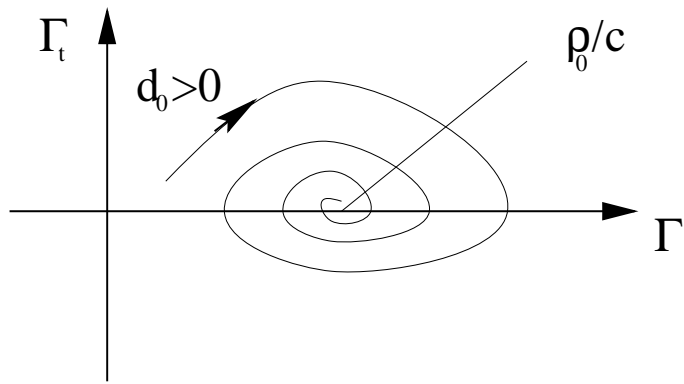

FIGURE 9. Phase plane graph of $(\Gamma,\Gamma_t)$ for $c > 0$ with damping

Thus, if we denote this first minimal time as $t^*(\alpha)$, then the global solution exists if and only if

$$\Gamma(\alpha,t^*) > 0.$$

In other words, the above expression of $\Gamma(\alpha, t)$ implies that there is a global in time solution if and only if,

$$\sqrt{(c-\rho_0)^2 + (\frac{cu_0'}{\mu} + \frac{c-\rho_0}{2\epsilon\mu})^2} < \rho_0 e^{\frac{t^*}{2\epsilon}},$$

and this is equivalent to the condition that

$$\left| u_0' + \frac{c-\rho_0}{2\epsilon c} \right| < \sqrt{\left(\frac{\mu\rho_0}{c}\right)^2 (e^{\frac{t^*}{\epsilon}} - 1) + \frac{\mu^2}{c}(2\rho_0 - c)}, \quad \forall \alpha \in \mathbb{R}.$$

Let us now compute $t^*$. By definition, $t^* = t^1$ is the smallest time such that

$$\Gamma_t(\alpha, t^n) = 0, \quad \Gamma_{tt}(\alpha, t^n) > 0, \quad n \in \mathbb{N},$$

where $0 \leq t^1 < t^2 < \cdots < t^n \to \infty$. ¿From $\Gamma_t(\alpha, t^n) = 0$ we have

$$\text{tg}(\mu t^n) = \frac{\mu \times \frac{2\epsilon c u_0' - \rho_0 + c}{2\epsilon\mu} - \frac{1}{2\epsilon}(c - \rho_0)}{(c-\rho_0)\mu + \frac{1}{2\epsilon} \times \frac{2\epsilon c u_0' - \rho_0 + c}{2\epsilon\mu}} = \frac{2\epsilon\mu u_0'}{u_0' + 2\epsilon k(c - \rho_0)}, \quad \mu^2 = ck - \frac{1}{4\epsilon^2}.$$

Substituting this $t^n$ into the expression of $\Gamma_{tt}$ we find

$$\Gamma_{tt}(\alpha, t^n) = \mu e^{-\frac{t}{2\epsilon}} \frac{\sin(\mu t^n)}{-u_0'} \left[ (u_0')^2 + \left( \frac{u_0'}{2\epsilon} + \frac{k}{\mu}(c - \rho_0) \right)^2 \right],$$

which shows that $t^n$ is a minimum with $\Gamma_{tt}(\alpha, t_n) > 0$ provided $\text{sgn}(u_0') = -\text{sgn}(\sin(\mu t^n))$. Therefore $t^*$ is uniquely determined by

(3.14)
$$t^* = \frac{1}{\mu} tg^{-1} \left[ \frac{2\epsilon\mu u_0'}{u_0' + 2\epsilon k(c-\rho_0)} \right] \begin{cases} 0 < t^* < \frac{\pi}{2\mu}, & \text{if } u_0' < \min\{0, 2\epsilon k(\rho_0 - c)\}, \\ \frac{\pi}{2\mu} < t^* < \frac{\pi}{\mu}, & \text{if } 2\epsilon k(\rho_0 - c) < u_0' < 0, \\ \frac{\pi}{\mu} < t^* < \frac{3\pi}{2\mu}, & \text{if } u_0' > \max\{0, 2\epsilon k(\rho_0 - c)\}, \\ \frac{3\pi}{2\mu} < t^* < \frac{2\pi}{\mu}, & \text{if } 0 < u_0' < 2\epsilon k(\rho_0 - c). \end{cases}$$

Moreover, we can now find the asymptotic behavior of $\rho(x,t)$ as $t \to \infty$. Indeed, following the same procedures as before, we find that $\rho(x(\alpha,t),t) = \frac{\rho_0(\alpha)}{\Gamma(\alpha,t)}$ with $\Gamma(\alpha, t)$ given in (3.14), and it follows that $\rho(x,t) \to c$ exponentially fast, as one would expect. We summarize by stating that

**Theorem 3.3.** *Consider the system of Euler-Poisson equations (3.1), (3.12), (3.3) with a constant background charge $c > 0$, a repulsive force $k > 0$ and weak relaxation $\epsilon > \frac{1}{2\sqrt{ck}}$. Let the critical time $t^*$ be defined in (3.14). Then if at all points $\alpha \in \mathbb{R}$,*

(3.15) $$\left| u_0'(\alpha) + \frac{c - \rho_0(\alpha)}{2\epsilon c} \right| < \sqrt{k - \frac{1}{4c\epsilon^2}} \sqrt{\frac{\rho_0^2(\alpha)}{c}(e^{t^*/\epsilon} - 1) + 2\rho_0(\alpha) - c},$$

*the solution of (3.1), (3.12), (3.3) is smooth for all time. In this case, $u(x,t) \to 0$ and $\rho(x,t) \to c$ exponentially as $t \to \infty$. Otherwise, if condition (3.15) fails, then the solution of (3.1), (3.12), (3.3) loses smoothness in a finite time.*

*Remark* 3.4. Note that as $\epsilon \to \infty$ we recover the local condition (3.10) for the case without relaxation.

Finally to complete our discussion of all choices of $k$, $c$ as well as $\epsilon$, we turn to consider the non-oscillatory case with strong relaxation $\epsilon < \frac{1}{2\sqrt{kc}}$. In this case we have

$$x'' = kE - \frac{x'}{\epsilon}.$$

Note that $E' = -cu = -cx'$, and hence
$$E = E_0(\alpha) - c(x - \alpha).$$
Combining the above two equations one finds that
$$x'' + \frac{1}{\epsilon}x' + ckx = k(E_0 + c\alpha),$$
and this equation together with the initial data, $(x, x')(\alpha, 0) = (\alpha, u_0(\alpha))$, leads to the characteristics of the form
$$x(\alpha, t) = \alpha + \frac{E_0}{c} + e^{-\frac{t}{2\epsilon}}\left\{[(-\frac{E_0}{c})(\frac{1}{2} - \frac{1}{4\epsilon\mu}) - \frac{u_0}{2\mu}]e^{-\mu t} + [(-\frac{E_0}{c})(\frac{1}{2} + \frac{1}{4\epsilon\mu}) + \frac{u_0}{2\mu}]e^{\mu t}\right\}.$$

Upon differentiation with respect to $\alpha$, one finds the indicator function $\Gamma(\alpha, t) = x_\alpha(\alpha, t)$,
$$\Gamma(\alpha, t) = \frac{\rho_0}{c} + e^{-\frac{t}{2\epsilon}}\left\{[(1 - \frac{\rho_0}{c})(\frac{1}{2} - \frac{1}{4\epsilon\mu}) - \frac{u_0'}{2\mu}]e^{-\mu t} + [(1 - \frac{\rho_0}{c})(\frac{1}{2} + \frac{1}{4\epsilon\mu}) + \frac{u_0'}{2\mu}]e^{\mu t}\right\}.$$

Based on this formula we claim the following.

**Theorem 3.4.** *Consider the system of Euler-Poisson equations (3.1), (3.12), (3.3) with a constant background charge $c > 0$, a repulsive force, $k > 0$ and a strong relaxation term, $\epsilon < \frac{1}{2\sqrt{ck}}$. If at all points $\alpha \in \mathbb{R}$*

$$(3.16) \qquad u_0'(\alpha) > \min\left\{0, -(1 - \frac{\rho_0(\alpha)}{c})\left(\sqrt{\frac{1}{4\epsilon^2} - ck} + \frac{1}{2\epsilon}\right)\right\},$$

*then the solution of (3.1), (3.12), (3.3) is smooth for all time. In this case, $u(x, t) \to 0$ and $\rho(x, t) \to c$ exponentially fast as $t \to \infty$.*

*Proof.* Expressed in terms of $\lambda = \frac{1}{2\epsilon}$, $a = \frac{1}{2}(1 - \frac{\rho_0}{c})$ and $b = (1 - \rho_0/c)\frac{1}{4\epsilon\nu} + \frac{u_0'}{2\nu}$ with $\nu = \sqrt{\frac{1}{4\epsilon^2} - ck}$, we have
$$\Gamma(\alpha, t) = 1 - 2a + (a - b)e^{-(\nu+\lambda)t} + (a + b)e^{(\nu-\lambda)t}.$$
Note that $\Gamma(\alpha, 0) = 1$ and $\Gamma(\alpha, \infty) = 1 - 2a = \frac{\rho_0}{c} > 0$. It suffices to show that $\Gamma(\alpha, t) > 0$ for all time $t > 0$ under the condition (3.16). First, if $a + b \geq 0$, then since $a < 1/2$ it follows that
$$\Gamma = (1 - 2a) + (a + b)e^{(\nu-\lambda)t}\left\{1 + \frac{a - b}{a + b}e^{-2\nu t}\right\}$$
$$\geq (1 - 2a) + (a + b)e^{(\nu-\lambda)t}\left\{1 + \min\left\{0, \frac{a - b}{a + b}\right\}\right\}$$
$$= (1 - 2a) + \min\{2a, a + b\}e^{(\nu-\lambda)t}$$
$$\geq \min\{1, 1 - 2a\} > 0.$$

Second, if $a + b < 0$, then one has to rule out a possibility of a local minimum achieved at positive time $t = t^* > 0$. At such time we would have $\Gamma_t(\alpha, t^*) = 0$, which gives
$$t^* = \frac{1}{2\nu}\log\left(\frac{(\nu + \lambda)(a - b)}{(\nu - \lambda)(a + b)}\right), \quad \Gamma_{tt}(\alpha, t^*) = 2\nu(a + b)(\nu - \lambda)e^{(\nu-\lambda)t^*} > 0.$$

Since $(\nu - \lambda)(a + b) > 0$, it follows that if $b > \frac{\lambda}{\nu}a$ then $(\nu + \lambda)(a - b)/(\nu - \lambda)(a + b) < 1$, and hence that $t^*$ is negative. Put differently, $t^* < 0$ and hence $\Gamma(\alpha, t) > 0$ for all $t > 0$, if $a + b < 0$ and $b > \frac{\lambda}{\nu}a$. In summary of these two cases, a sufficient condition for the global existence is $b > \min\{\frac{\lambda}{\nu}a, -a\}$, which is exactly the same as (3.16) when recalling the definition of $a$ and $b$. □

## 4. CRITICAL THRESHOLDS– 1D MODEL WITH VISCOSITY

We consider the solutions of a parabolic- hyperbolic version of (2.1)-(2.3)–the modified viscous Burgers-Poisson equations:

(4.1) $$\rho_t + (\rho u)_x = 0,$$

(4.2) $$u_t + uu_x = -k\phi_x + \left(\frac{u_x}{\rho}\right)_x = kE + \left(\frac{u_x}{\rho}\right)_x,$$

(4.3) $$E_x = -\phi_{xx} = \rho,$$

with repulsive force $k > 0$ (if the matter is treated as charged particles). We show that despite the presence of a parabolic term on the right of (4.2), these equations can still lose smoothness in finite time when the critical threshold condition is crossed.

Note that the parabolic term on the right of (4.2), $\left(\frac{u_x}{\rho}\right)_x$ is similar to the one, $u_{xx}/\rho$, found in the Navier-Stokes equations, $u_t + uu_x = u_{xx}/\rho$, with a difference term of $-\frac{\rho_x u_x}{\rho^2}$. In this context, both types of viscosity regularization terms admit similar behavior. As $\rho \to 0$ then $u_{xx}/\rho$, $-\frac{\rho_x u_x}{\rho^2}$ and hence $\left(\frac{u_x}{\rho}\right)_x$ all tend to infinity. More importantly, the forcing effect in (4.2) leads to $u \downarrow -\infty$, which in turn leads by (4.1), that $\rho \to \infty$, and then $u_{xx}/\rho$, $-\frac{\rho_x u_x}{\rho^2}$ and hence $\left(\frac{u_x}{\rho}\right)_x$ all tend to zero. Thus viscous terms tend to zero when they are most needed and the solution blows up with $\rho \uparrow \infty$. It is this vanishing viscosity mechanism which allows the possible blow up above the critical threshold in Euler-Poisson system (4.1)-(4.3).

Let $d = u_x$. Differentiate (4.2) with respect to $x$ to obtain

(4.4) $$d_t + ud_x = -d^2 + k\rho + (d/\rho)_{xx}.$$

The difference, $\rho \times (4.4) - d \times (4.1)$ yields

(4.5) $$\rho d_t - d\rho_t + u\rho d_x - u\rho_x d = k\rho^2 + \rho(d/\rho)_{xx}.$$

Our goal as before is to estimate the ratio $\beta := d/\rho$, thus $\frac{1}{\rho^2} \times (4.5)$ yields the following parabolic equation

$$\beta_t + u\beta_x = k + \beta_{xx}/\rho.$$

We would like to show a $\beta$-bound due to a maximum principle of the form,

(4.6) $$\inf_\alpha \beta(\alpha, 0) + kt \leq \beta(x, t) \leq \sup_\alpha \beta(\alpha, 0) + kt.$$

The difficulty here is that $\rho(x, t)$ must die at infinity if one wants the integral of $d$, that is if one wants $u$, to be finite at $\pm\infty$ and one wants $\beta$ to be bounded away from 0. If $\rho$ dies at infinity, then viscosity has a coefficient that blows-up at infinity. We make use of the following theorem to show that even in such cases we still have a maximum principle.

**Theorem 4.1** (maximum principle). *If $u_t + f(x,t)u_x = a(x,t)u_{xx}$, $u(x,t) \leq D(x)E(t)$ where $D(x)$ is sub-linear in $x$ and $E(t)$ is exponential in time; $u, u_t, u_x, u_{xx} \in C$, $|f(x,t)| \leq d$; and $0 \leq a(x,t)$; $a(x,t) \in C$, then the solution of the equation obeys the maximum principle:*

$$u(x,t) \leq \sup_x u(x,0).$$

To the best of our knowledge there seems no maximum principle of this form available in the literature. We present a self-contained proof of this theorem in the Appendix.

To proceed we assume at most sub-linear growth of $u_x/\rho$ at far field in space and exponential growth of $u_x/\rho$ in time; i.e. there exists $\gamma > 0$ such that for $|x| \gg 1$

$$\left|\frac{u_x}{\rho}\right| \ll \text{Const.} e^{\gamma t}|x|. \tag{4.7}$$

Equipped with the above maximum principle we establish the following

**Theorem 4.2.** *Consider the system (4.1)- (4.3) with the repulsive force $k > 0$. Assume that the smooth solution is sought such that $u_x/\rho$ is sub-linear in the sense of satisfying (4.7). Then if for any $\alpha \in \mathbb{R}$*

$$\sup_x \left(\frac{u'_0(x)}{\rho_0(x)}\right) < -\sqrt{\frac{2k}{\rho_0(\alpha)}},$$

*then $(u, \rho)$ ceases to behave nicely in finite time. If*

$$\inf_x \left(\frac{u'_0(x)}{\rho_0(x)}\right) > -\sqrt{\frac{2k}{\rho_0(\alpha)}},$$

*then $(u, \rho)$ remains smooth for all time. In this case as $t \to \infty$*

$$\rho(x,t) \sim \frac{2}{kt^2}, \quad d(x,t) \sim \frac{2}{t}.$$

*Proof.* Theorem 4.1 shows that $\beta$ satisfies a maximum principle as long as $u \in C^3$, $|u|$ is bounded, $\rho \in C^2$, and $u_x/\rho$ is sub-linear in $x$ and exponential in $t$. We assumed that these conditions all hold, and we make use of the maximum principle to bound $\beta$ by (4.6).

Because (4.1) says that the derivative of $\rho$ along the curve $x(\alpha, t)$ is $-\rho d$, that is,

$$\rho' = -\rho d = -\beta(t)\rho^2,$$

where $\beta(t) = d/\rho$ satisfies the bounds stated in (4.6). By Lemma 2.1 one has

$$\rho(x(\alpha,t),t) = \frac{\rho_0(\alpha)}{1 + \rho_0(\alpha)\int_0^t \beta(\tau)d\tau}.$$

This gives

$$\frac{\rho_0(\alpha)}{\rho_0(\alpha)(kt^2/2 + \sup_x \beta(x,0)t) + 1} \leq \rho(x(\alpha,t),t) \leq \frac{\rho_0(\alpha)}{\rho_0(\alpha)(kt^2/2 + \inf_x \beta(x,0)t) + 1}.$$

Finally, as we have bounds for $\beta = d/\rho$ we find that $d$ must satisfy:

$$\frac{\rho_0(\alpha)(kt + \inf_x \beta(x,0))}{\rho_0(\alpha)(kt^2/2 + \inf_x \beta(x,0)t) + 1} \leq d(x(\alpha,t),t) \leq \frac{\rho_0(\alpha)(kt + \sup_x \beta(x,0))}{\rho_0(\alpha)(kt^2/2 + \sup_x \beta(x,0)t) + 1}$$

as long as both of the numerators remain positive. Suppose that $k = 1$, the initial data $\rho_0 = \frac{1}{1+x^2}$, $u_0 = -2\arctan(x)$. (Note that the total mass is finite.) We find that $\beta(x,0) \equiv -2$. Our maximum principle argument says that either $d \to -\infty$ as $t \to 2 - \sqrt{2}$ or $(u, \rho)$ ceases to behave well sometime before $t = 2 - \sqrt{2}$. The "weakest" way in which $\beta$ can cease to behave well is for $\beta = u_x/\rho$ to grow faster than sub-linearly in space. As $\beta(x,0) \equiv -2$, we find the supremum and infimum of $\beta$ are equal. The above inequalities lead to equalities until the maximum principle no longer applies. We find that (as $k = 1$):

$$\rho(x(\alpha,t),t) = \frac{\rho(\alpha)}{\rho(\alpha)(t^2/2 - 2t) + 1}.$$

As long as $\rho(x,t)$ is continuous in time this implies that up to and including the time $t_1 < 2 - \sqrt{2}$ at which $\beta$ stops being sub-linear in space $\rho(x(\alpha,t),t)$ is completely known. Moreover, for large $\alpha$ we know that $\rho(x(\alpha,0),0) \approx 1/\alpha^2$. (As by assumption $|u(x,t)|$ is bounded, we know that $x(\alpha,t)$ does not change quickly.) Thus if $\alpha$ is large enough we find that $\rho(x(\alpha,t),t) \approx \rho(\alpha) \approx 1/\alpha^2$. This is sufficient to guarantee that $\rho(x,t)$ does not decay faster than $1/x^2$. In other words, if $\rho(x,t)$ decays like $1/x^2$, it cannot decay any faster at a later time. If $u_x/\rho$ grows non sub-linearly, then $u_x$ will die at a rate that is slower than $1/(\ln(x)x)$—that is $u$ will tend to infinity.

We note that if we take initial data in which $\rho_0(\alpha)(kt^2/2 + \sup_x \beta(x,0)t) + 1$ is ever zero–if we ever cross this critical threshold–then the solution must either blow up in finite-time (because the solution has an upper bound that tends to $-\infty$ in finite- time) or the solution must cease to be nicely behaved in finite-time (so that the maximum principle ceases to apply). In either case, we find that once the critical threshold is crossed, the solution no longer behaves nicely for all time.

The estimates above show that as long as $k > 0$ and the denominators do not ever become zero, $d \to 2/t$ and $\rho \to 2/(kt^2)$. We see that the solution of the parabolic-hyperbolic equation tends to behave like the solution of the hyperbolic one in the long term. Similarly, we find that if $k = 0$ and the denominators remain positive for all time, then:

$$\frac{c_1}{t} \leq d \leq \frac{c_2}{t},$$

where $c_1 \leq 1 \leq c_2$. This is also the behavior that one would expect without the parabolic term (i.e. the inviscid Burgers equation). Hence the proof is complete. □

We see that the presence of a self-induced electric field and of viscosity is not necessarily enough to stop singularities from forming. This is not terribly surprising. The viscosity used tends to zero as $\rho \to \infty$–it is smallest when it needs to be largest. This, of course, is a problem with the viscous term in the Navier-Stokes equations as well. Additionally, as we are dealing with equations in one dimension, the electric field that we are dealing with is essentially the field due to an infinite sheet of charge. The field due to such a sheet is a constant throughout space–forcing two such sheets to collide is easy. There is no reason that such a force should prevent the collision of sheets. What we have succeeded in doing is showing how the behavior of the solutions of our sets of equations reflects the behavior of the physical systems they represent.

## 5. CRITICAL THRESHOLDS–MULTI-D MODEL WITH GEOMETRICAL SYMMETRY

Let us consider the Euler-Poisson equations governing the $\nu + 1$ dimensional isotropic ion expansion in the electrostatic fluid approximation for cold ions

(5.1) $$r^\nu n_t + (nur^\nu)_r = 0,$$

(5.2) $$u_t + uu_r = -k\phi_r = kE,$$

(5.3) $$(r^\nu \phi_r)_r = -nr^\nu,$$

subject to the initial data

(5.4) $$(n,u)(r,0) = (n_0, u_0)(r), \quad n_0(r) \geq 0.$$

Here $r > 0$ denotes the distance from the origin, and $k > 0$ is a known physical constant, which takes into account the general scaled version (5.3) of Poisson equation $(r^\nu \phi_r)_r =$ $-4\pi q n r^\nu$. The unknowns are the local particle density $n = n(r,t)$, velocity $u = u(r,t)$ and potential $\phi = \phi(r,t)$. The geometrical factor $\nu$ takes values $0, 1, 2$, for planar, cylindrical, or spherical symmetry.

Our interest is to find a critical threshold criterion for the Cauchy problem (5.1)-(5.4). Several issues that will also be clarified along the way include:
(1) Expansion rate of the flow path;
(2) The large time behavior of the velocity;
(3) The decay rate of the density as well as the velocity gradient;
(4) Sharp estimate of blow up time when the initial data exceed the critical threshold.

### 5.1. Analytic solution along the particle path.
Without pressure forces, the particles move without collisions in a self-consistent electrostatic field. The system contains only particle-path characteristics. We shall trace the time dynamics along these characteristics.

To this end we introduce the following weighted quantities

$$e := Er^\nu = -\phi_r r^\nu, \quad \rho = nr^\nu.$$

With these definitions, (5.1) becomes

(5.5) $$\rho_t + (\rho u)_r = 0,$$

and the Poisson equation (5.2), now reads

(5.6) $$e_r = \rho.$$

To avoid having a singularity at the origin we require $e(0,t) = 0$, hence

$$e = \int_0^r \rho(\xi, t) d\xi.$$

In view of (5.5), this implies that $e$ satisfies the transport equation

(5.7) $$e_t + u e_r = 0.$$

To solve this equation we define the "flow map" $r(\alpha, t) : \mathbb{R}^+ \to \mathbb{R}^+$

(5.8) $$\frac{dr(\alpha, t)}{dt} = u(r(\alpha, t), t), \quad r(\alpha, 0) = \alpha.$$

Denoting differentiation along this characteristic curve by $' := \frac{d}{dt}$, the mass equation (5.7) and the momentum equation (5.2) yield

(5.9) $$e' = 0,$$

(5.10) $$u' = \frac{ke}{r^\nu}, \quad r' := u(r, t).$$

We now solve the above system subject to the initial data

$$(r, e, u)|_{t=0} = (\alpha, e_0(\alpha), u_0(\alpha))$$

with $\alpha \in \mathbb{R}^+$ parameterizing the initial location, and with the weighted mass $e_0(\alpha) = \int_0^r n_0(\xi) \xi^\nu d\xi$.

By (5.9), $e$ remains constant along the isotropic characteristics, $e = e_0(\alpha)$, and (5.10) then yields

(5.11) $$r'' = \frac{ke_0}{r^\nu}, \quad r(\alpha, 0) = \alpha, \quad r'(\alpha, 0) = u_0(\alpha).$$

This relation shows that each fluid particle, starting from position $\alpha$ with initial velocity $u_0(\alpha)$, is influenced by a central acceleration $\frac{ke_0}{r^\nu}$. Once the solution of this equation for $r(\alpha, t)$ is known, the other dependent variables $u, \rho$ can be determined accordingly. We proceed to study the solution of (5.11). To this end, we introduce the 'indicator' function

$$\Gamma(\alpha, t) := e^{\int_0^t u_r(r(\alpha,\tau),\tau)d\tau}.$$

The geometrical interpretation of $\Gamma(\alpha, t)$ will be clear from the explicit solution of the Euler- Poisson system(5.1)-(5.3), given in the following lemma which will play an essential role in our discussion.

**Lemma 5.1.** *Consider the Euler-Poisson equations (5.1)- (5.3) subject to the initial data $(n_0, u_0) \in C^1(\mathbb{R}^+) \times C^1(\mathbb{R}^+)$. Let $r(\alpha, t)$ be the flow map defined in (5.8), then*

$$\Gamma(\alpha, t) = \frac{\partial r(\alpha, t)}{\partial \alpha}.$$

*Moreover, the solution of (5.1)- (5.4) is given by*

$$(5.12) \qquad u(r, t) = \frac{\partial r(\alpha, t)}{\partial t},$$

$$(5.13) \qquad n(r, t) = \frac{n_0(\alpha)\alpha^\nu}{r^\nu \Gamma(\alpha, t)},$$

$$(5.14) \qquad u_r(r, t) = \frac{\Gamma_t(\alpha, t)}{\Gamma(\alpha, t)}.$$

*Proof.* Along the particle path one has

$$\frac{d}{dt}r(\alpha, t) = u(r(\alpha, t), t), \quad r(\alpha, 0) = \alpha, \quad \forall \alpha \in \mathbb{R}^+.$$

Differentiating this equation with respect to $\alpha$ gives

$$\frac{d}{dt}\left(\frac{\partial}{\partial \alpha}r(\alpha, t)\right) = u_r(r(\alpha, t), t)\frac{\partial}{\partial \alpha}r(\alpha, t), \quad \frac{\partial}{\partial \alpha}r(\alpha, 0) = 1.$$

Hence $\frac{\partial}{\partial \alpha}r(\alpha, t) = \Gamma(\alpha, t)$ for any $t \in \mathbb{R}^+$.

The mass equation along the particle path $r(\alpha, t)$ becomes

$$\rho' + \rho u_r = 0,$$

and integration in time leads to

$$\rho(r, t) = \rho_0(\alpha)/\Gamma(\alpha, t).$$

The formula for $u_r$ follows from the definition of $\Gamma(\alpha, t)$. □

Using the expression of the solution, (5.12)- (5.14), we conclude

**Corollary 5.2.** *The smooth solution to the Euler- Poisson equations (5.1)-(5.4) blows up in finite time, $t = t_c$, if and only if one of the following equivalent conditions is met:*
*(1) $\int_0^{t_c} u_r(r(\alpha, \tau), \tau)d\tau = -\infty$;*
*(2) $\Gamma(\alpha, t_c) = 0$;*
*(3) There exists an $\alpha \in \mathbb{R}$ such that $\frac{\partial r}{\partial \alpha}(\alpha, t_c) = 0$.*

To ensure the existence of the global regular solution, therefore it suffices to start with the set of prescribed initial data for which ( recall $\frac{\partial r}{\partial \alpha}(\alpha, 0) = 1$),

$$\frac{\partial r}{\partial \alpha}(\alpha, t) > 0 \quad \forall t \in \mathbb{R}^+.$$

5.2. **The isotropic flow map.** Equipped with the above relations we are in a position to study the isotropic flow map $r(\alpha, t)$, and the zeros of $r_\alpha(\alpha, t)$ which characterize the formation of the singularity.

We begin with the solution of the isotropic flow map $r(\alpha, t)$ governed by (5.11). We summarize its behavior in the following Lemma.

**Lemma 5.3.** *The solution of d-dimensional problem, $r'' = ke_0(\alpha)r^{-\nu}$, with initial data $(r, r')(\alpha, 0) = (\alpha, u_0(\alpha) > 0)$ is as follows. (We classify the different cases by the value of $\nu := d - 1$.*

- $\nu = 0$: *the flow map is given by*

(5.15) $$r(\alpha, t) = \alpha + u_0(\alpha)t + \frac{ke_0(\alpha)}{2}t^2,$$

*and the velocity is*

(5.16) $$u(r, t) = u_0(\alpha) + ke_0(\alpha)t.$$

- $\nu = 1$: *the flow map is*

(5.17) $$r(r, t) = \alpha \exp\left(\frac{ke_0}{2}\tau^2 + u_0\tau\right),$$

*with the parameter $\tau$ given implicitly by the temporal relation*

$$t = \alpha \int_0^\tau \exp\left(\frac{ke_0}{2}\xi^2 + u_0\tau\right) d\xi.$$

*In this case we have*

(5.18) $$\sqrt{\alpha^2 + 2\alpha u_0 t + ke_0 t^2} \leq r(\alpha, t) \leq \alpha + u_0 t + O(t\ln t).$$

*The corresponding velocity is given by*

(5.19) $$u(r, t) = h(\alpha, t),$$

*with $h(\alpha, t) \sim u_0 + O(t\ln t)$ determined implicitly by the identity*

$$t \equiv \int_{u_0(\alpha)}^{h(\alpha,t)} \exp\left(\frac{\xi^2 - u_0^2}{2ke_0}\right) d\xi.$$

- $\nu = 2$: $r = r(\alpha, t)$ *is given implicitly by*

(5.20) $$\frac{2Q(\alpha)}{R}t = \cosh^{-1}(\frac{2r}{R} - 1) - \cosh^{-1}(\frac{u_0^2}{ke_0}) + \frac{2}{R}\sqrt{r^2 - Rr} - \frac{u_0}{\sqrt{ke_0}}.$$

*with*

(5.21) $$Q(\alpha) := \sqrt{u_0^2 + \frac{2ke_0}{\alpha}}, \quad R(\alpha) := \frac{2ke_0}{u_0^2 + \frac{2ke_0}{\alpha}}.$$

*In this case we have*

(5.22) $$[\alpha^3 + 3\alpha^2 u_0 t + \frac{3}{2}ke_0 t^2]^{1/3} \leq r(\alpha, t) \leq \alpha + \sqrt{u_0^2 + \frac{2ke_0}{\alpha}}\, t,$$

*the velocity is uniformly bounded in time, and*

(5.23) $$\lim_{t \to \infty} u(r, t) = Q(\alpha).$$

*Proof.* We start with $\nu = 0$. A straightforward integration of $r'' = ke_0(\alpha)$ combined with the initial position $\alpha$ and the initial velocity $u_0(\alpha)$ gives the same formula for $r$ that we met earlier, (2.10), in the 1- dimensional case.

We then turn to the 2-dimensional case, $\nu = 1$, where the flow map equation (5.11) reads

$$(5.24) \qquad r'' = ke_0(\alpha)r^{-1}.$$

Starting at the location $\alpha$ with velocity $u_0$, then the energy integral is

$$(5.25) \qquad \frac{1}{2}(r')^2 - \frac{1}{2}u_0^2 = ke_0(\ln r - \ln \alpha),$$

implying

$$(5.26) \qquad r' = \sqrt{u_0^2 + 2ke_0 \ln(r/\alpha)}, \quad u_0 > 0.$$

Let $\tau(\alpha, t)$ be a dimensionless parameter such that for constant $\alpha$

$$d\tau = \frac{dt}{r} = \frac{dr}{r\sqrt{u_0^2 + 2ke_0 ln(r/\alpha)}} = \frac{1}{ke_0}d\left[\sqrt{u_0^2 + 2ke_0 ln(r/\alpha)}\right], \quad r(\tau = 0) = \alpha.$$

Then the parametric solution of (5.26) is given by (5.17), i.e.,

$$r = \alpha \exp\left(\frac{ke_0}{2}\tau^2 + u_0\tau\right).$$

Here, the parameter $\tau$ is determined by $t$

$$dt = rd\tau = \alpha \exp\left(\frac{ke_0}{2}\tau^2 + u_0\tau\right)d\tau,$$

which gives

$$t = \alpha \int_0^\tau \exp\left(\frac{ke_0}{2}\xi^2 + u_0\xi\right)d\xi.$$

We want to show that $r \sim O(t)$ for large $t$, with the tight bound given by (5.18). To show the lower bound, we rewrite the equation (5.24) as

$$(rr')' = ke_0 + r'^2 \geq ke_0,$$

and integration twice gives the bound on the left of (5.18). Combining this lower bound with $r'' = ke_0 r^{-1}$ yields

$$r'' \leq ke_0[\alpha^2 + 2\alpha u_0 t + ke_0 t^2]^{-1/2},$$

and integration twice gives the upper bound shown on the right of (5.18).

Finally we turn to the 3-dimensional case $\nu = 2$, where the flow map equation reads

$$(5.27) \qquad r'' = ke_0(\alpha)r^{-2}.$$

Starting at location $\alpha$ with velocity $u_0$, the energy integral is

$$(5.28) \qquad \frac{1}{2}r'^2 - \frac{1}{2}u_0^2 = ke_0(\frac{1}{\alpha} - \frac{1}{r}),$$

which gives

$$r' = Q(\alpha)(1 - R/r)^{1/2} \quad \text{for} \quad u_0 \geq 0,$$

with $Q$ and $R$ given in (5.21). Integration in time once gives the implicit formula of the flow map (5.22), where the dependence of $r$ on $\alpha$ will play essential roles (i.e., $r_\alpha = 0$) in our later analysis. We conclude with the large time estimate (5.22). The upper-bound follows from the fact $r' \leq Q(\alpha)$. ¿From the flow map equation (5.27) follows that
$$(r^3/3)'' = 2r(r')^2 + ke_0 \geq ke_0,$$
and integration twice yields the bound on the left of (5.22). □

*Remark* 5.1. The above formula reveals quite different geometry of the isotropic flow path. In the case $\nu \geq 1$, the energy integral ensures the lower positive bound for $r$. In fact recalling (5.25) and (5.28),

$$(5.29) \qquad (r')^2 = u_0^2 + 2ke_0 \ln r/\alpha, \quad \nu = 1,$$

$$(5.30) \qquad (r')^2 = u_0^2 + \frac{2ke_0}{\nu - 1}(\alpha^{1-\nu} - r^{1-\nu}), \quad \nu > 1,$$

one finds

$$(5.31) \qquad r(\alpha, t) \geq \alpha e^{-\frac{u_0^2}{2ke_0}}, \quad \nu = 1$$

$$(5.32) \qquad r(\alpha, t) \geq \left[\frac{\nu - 1}{2ke_0}u_0^2 + \alpha^{1-\nu}\right]^{-\frac{1}{\nu-1}}, \quad \nu > 1.$$

Moreover following the same procedure as for $\nu = 2$ we obtain the expansion rate of the flow path for general situation $\nu \geq 2$,

$$(5.33) \quad [\alpha^{\nu+1} + (\nu+1)\alpha^\nu u_0 t + \frac{\nu+1}{2}ke_0 t^2]^{1/(\nu+1)} \leq r(\alpha, t) \leq \alpha + \sqrt{u_0^2 + \frac{2ke_0}{(\nu-1)\alpha^{\nu-1}}}\, t.$$

In the 1D case, $\nu = 0$, we saw that if the initial velocity is negative then particles can reach the $r = 0$ line in a finite time due to the quadratic form (5.15). In contrast, for $\nu \geq 1$, these positive lower bounds imply that if the velocity is initially negative, the particle path $r$ may decrease only for a finite time and then increase due to the positive acceleration. To avoid the technical discussions here and in what follows we restrict ourselves to the case $u_0 > 0$.

*Remark* 5.2. The above results show that for the case $\nu = 0, 1$ the velocity grows linearly. But for the case $\nu = 2$, the velocity is uniformly bounded and converges to a positive constant as time becomes large.

To study the critical threshold phenomena one may utilize the 'indicator' function $\Gamma = r_\alpha$. This study is carried out in next section.

**5.3. Critical thresholds.** Using the above flow map we provide precise conditions on the initial data such that either the solution remains globally smooth or it breaks down in a finite time.

We start by revisiting the 1-dimensional case.

**Theorem 5.4** (Global existence of smooth solutions for planar case ($\nu = 0$)). *The smooth solutions of (5.1)-(5.4) with $\nu = 0$ exist if and only if*

$$(5.34) \qquad u_0'(\alpha) > -\sqrt{2kn_0(\alpha)}, \quad \forall \alpha \in \mathbb{R}^+.$$

*In this case the solution is given by*

$$n(r(\alpha,t),t) = \frac{n_0(\alpha)}{1 + u_0't + \frac{k}{2}n_0(\alpha)t^2},$$

$$u_r(r(\alpha,t),t) = \frac{u_0' + kn_0 t}{1 + u_0't + \frac{k}{2}n_0(\alpha)t^2}.$$

*Proof.* Differentiating the flow map equation $r = ke_0$ with respect to $\alpha$ we find that

$$\Gamma'' = kn_0(\alpha),$$

for $e_{0\alpha} = \rho_0 = n_0$. The definition of $\Gamma$ gives $(\Gamma', \Gamma)(t=0) = (u_0', 1)$. Thus the corresponding energy integral is

$$[\Gamma']^2 - 2kn_0\Gamma = [u_0']^2 - 2kn_0.$$

The geometry of the trajectory implies that to ensures the positivity of $\Gamma$ the initial data should satisfy either $u_0' \geq 0$ or $[u_0']^2 - 2kn_0 < 0$ for the case $u_0' < 0$, which yields (5.34). □

The above formula immediately yields the following

**Corollary 5.5** (Breakdown of smooth solutions for planar case ($\nu = 0$)). *The smooth solution to the Euler-Poisson equations (5.1)-(5.3) blows up in finite time if and only if the condition, $u_0'(\alpha) > -\sqrt{2kn_0(\alpha)}$, fails, i.e., if*

$$\exists \alpha \in \mathbb{R}^+ \quad s.t. \quad u_0'(\alpha) \leq -\sqrt{2kn_0(\alpha)}.$$

*In this case, the density $n(r,t)$ and $u_r(r,t)$ become infinite as $t \uparrow T$, where the blow-up time $t = t_c$, is given explicitly by*

$$t_c := 2/\sup\left\{-u_0' + \sqrt{(u_0')^2 - 2kn_0}\right\}.$$

*Remark* 5.3. Consider the 1D equation with nonzero background and additional relaxation term, where the equation (5.2) is replaced by

$$u_t + uu_r = kE - u/\epsilon \quad \text{with} \quad E_r = (n - c).$$

Hence the equation for $r$ reads

$$r'' = kE - \frac{r'}{\epsilon}.$$

Note that $E' = -cu = -cr'$, thereby $E = E_0(\alpha) - c(r - \alpha)$. These equations lead to an 'indicator' function $\Gamma = r_\alpha$, satisfying

$$\Gamma'' + \frac{\Gamma'}{\epsilon} + kc\Gamma = k\rho_0.$$

Using the phase plane analysis one recovers, for the 1D half space problem, $\alpha \in \mathbb{R}^+$, the same results obtained for the 1-D Cauchy problem, $\alpha \in \mathbb{R}$, consult Section 2-3.

The rest of this section is devoted to the multi- D case, $\nu \geq 1$, where we confirm the remarkable presistence of critical threshold phenomena in the multidimensional problem. But precisely determining sharp critical thresholds is far from trivial. Even in the isotropic case, the implicit solution formula makes the final conditions on initial data rather cumbersome. Thus, for example, differentiating (5.11) with respect to $\alpha$ yields

$$\Gamma'' = k\rho_0 r^{-\nu} - k\nu e_0 r^{-(\nu+1)}\Gamma,$$

which is coupled with the flow map equation $r'' = ke_0 r^{-\nu}$. It is difficult to find an explicit sharp threshold for the initial data that distinguishes between cases for which $\Gamma$ remains nonzero and cases for which it does not. Though some further tedious calculations may enable us to obtain a complex criterion for cylindrical case ($\nu = 1$) as well as the spherical case ($\nu = 2$), we do not perform these calculations. Instead we give sufficient conditions for upper thresholds on the initial data for the existence of global smooth solution, as well as the lower thresholds for the finite time breakdown. These confirm the existence of an intermediate critical threshold —which is the focus of our interest in this work.

We start with the 2D case.

**Theorem 5.6** (Global existence of smooth solutions for the cylindrical case $\nu = 1$). *A global smooth solution of Euler- Poisson equations (5.1)-(5.3) with $\nu = 1$ exists provided the initial data $(u_0, n_0)$ with $E_0 = \alpha^{-1} \int_0^\alpha n_0(\xi) \xi d\xi$ satisfy*

$$(5.35) \qquad u_0' > -\frac{k}{u_0}[\alpha n_0 h(\alpha) - E_0], \quad \forall \alpha \in \mathbb{R}^+,$$

*where $h(\alpha)$ is determined by*

$$(5.36) \qquad k\alpha^2 n_0 u_0 \int_0^{h(\alpha)} \frac{h(\alpha) - \eta}{[u_0^2 + 2kE_0 \alpha \eta]^{3/2}} e^\eta d\eta \equiv 1, \quad \forall \alpha \in \mathbb{R}^+.$$

*Proof.* Recall the energy integral (5.26)

$$r' = \sqrt{u_0^2 + 2ke_0 ln(r/\alpha)},$$

from which it follows that

$$\int_\alpha^{r(\alpha,t)} \frac{d\xi}{\sqrt{u_0^2 + 2ke_0 ln(\xi/\alpha)}} = t.$$

Differentiating the above equality with respect to $\alpha$, one has

$$\frac{\Gamma(\alpha,t)}{\sqrt{u_0^2 + 2ke_0 ln(r/\alpha)}} = A(\alpha,t),$$

where

$$A(\alpha,t) := \int_\alpha^{r(\alpha,t)} \frac{u_0 u_0' + k\rho_0 ln\xi/\alpha - ke_0 \alpha^{-1}}{[u_0^2 + 2ke_0 ln\xi/\alpha]^{3/2}} d\xi + \frac{1}{u_0}.$$

Noting that $e_0 = \alpha E_0$, we introduce $M := [kE_0 - u_0 u_0']/(k\alpha n_0)$ and rewrite $A$ in terms of $\eta := ln(\xi/\alpha)$, as

$$A(\alpha,t) = k\alpha^2 n_0 \int_0^{ln(r/\alpha)} \frac{\eta - M}{[u_0^2 + 2kE_0 \alpha \eta]^{3/2}} e^\eta d\eta + \frac{1}{u_0}.$$

We shall show that $A$, and therefore $\Gamma$, remain positive for all $t > 0$ provided (5.35) holds. W consider the case $M > 0$, since the complementary case $M \leq 0$ is trivial. Note that $e_{0\alpha} = \rho_0 = n_0 \alpha$. A simple computation involving (5.26) gives

$$\frac{dA}{dt} = k\alpha n_0 \frac{ln(r/\alpha) - M}{u_0^2 + 2kE_0 \alpha ln(r/\alpha)}.$$

In view of the monotonicity of $r$ in time ($\frac{dr}{dt} > 0$) we see that $A$ may achieve its unique minimum at $t^*$, where $A_t(t^*) = 0$, i.e., $r(\alpha, t^*) = \alpha e^M$, and this minimum is positive provided (5.35) holds, for

$$A(\alpha, t^*) = k\alpha^2 n_0 \int_0^M \frac{\eta - M}{[u_0^2 + 2kE_0\alpha\eta]^{3/2}} e^\eta d\eta + \frac{1}{u_0}$$

$$> k\alpha^2 n_0 \int_0^{h(\alpha)} \frac{\eta - h(\alpha)}{[u_0^2 + 2kE_0\alpha\eta]^{3/2}} e^\eta d\eta + \frac{1}{u_0} = 0,$$

if $M < h(\alpha)$, which is equivalent to (5.35). Therefore the indicator function $\Gamma(\alpha, t) = A(\alpha, t)\sqrt{u_0^2 + 2ke_0 ln(r/\alpha)}$ remains positive because $A(\alpha, t) \geq A(\alpha, t^*) > 0$ for all $t > 0$. □

Condition (5.35) could be viewed as an upper threshold in the sense of providing a sufficient condition leading to global smooth solutions, though the permissible class of the initial data for global smooth solutions is clearly larger. However, the existence of the critical threshold can be ensured by combining this upper threshold with the following lower threshold for the finite time breakdown.

**Theorem 5.7** (Breakdown of smooth solutions for cylindrical case $\nu = 1$). *The smooth solution to the Euler- Poisson equations (5.1)-(5.3) with $\nu = 1$ breaks down in finite time if the condition, $u_0'(\alpha) > -\sqrt{2kn_0(\alpha)}$, fails, i.e., if*

(5.37) $$\exists \alpha \in \mathbb{R}^+ \quad s.t. \quad u_0'(\alpha) \leq -\sqrt{2kn_0(\alpha)}.$$

*Proof.* Using the parametric form of the flow map given in (5.17) we evaluate the 'indicator' function $\Gamma$ via

$$\Gamma(\alpha, t) = r_\alpha - r_\tau/t_\tau t_\alpha.$$

¿From (5.17) we see that

$$r_\alpha = \left(1 + \alpha(\frac{k}{2}\rho_0\tau^2 + u_0'\tau)\right) \exp\left[\frac{ke_0}{2}\tau^2 + u_0\tau\right],$$

$$r_\tau = \alpha(ke_0\tau + u_0) \exp\left[\frac{ke_0}{2}\tau^2 + u_0\tau\right],$$

$$t_\alpha = \int_0^\tau \left(1 + \alpha(\frac{k}{2}\rho_0\xi^2 + u_0'\xi)\right) \exp\left[\frac{ke_0}{2}\xi^2 + u_0\xi\right] d\xi,$$

$$t_\tau = \alpha\exp\left(\frac{ke_0}{2}\tau^2 + u_0\tau\right) = r.$$

Expressed in terms of

(5.38) $$b(\tau) := 1 + \alpha u_0'\tau + \frac{k}{2}\alpha\rho_0\tau^2$$

and $r = r(\tau) = \alpha\exp[\frac{ke_0}{2}\tau^2 + u_0\tau]$, the 'indicator' function can be rewritten as

$$\Gamma = \frac{r}{\alpha}\left\{b(\tau) - (u_0 + ke_0\tau)\frac{\int_0^\tau r(\xi)b(\xi)d\xi}{r(\tau)}\right\}.$$

Note that $b(0) = 1$. The quadratic form of (5.38) implies that there must exist a parameter $\tau^*$ such that $b(\tau^*) = 0$ provided (5.37) holds. At this time $\Gamma$ becomes negative because the nonlocal term $\int_0^{\tau^*} r(\xi)b(\xi)d\xi$ stays positive. This combined with the fact $\Gamma(0) = 1$ ensures that there must be a finite time $t = t^*$ such that $\Gamma(\alpha, t^*) < 0$. Hence $\Gamma$ must vanish at finite time $t = t_c < t^*$. This completes the proof. □

We conclude with the 3-dimensional case, stating the lower threshold for finite time breakdown.

**Theorem 5.8** (Breakdown of smooth solutions for spherical case $\nu = 2$). *The solution of Euler-Poisson equations (5.1)-(5.3) for $\nu = 2$ blows up in finite time if the condition, $u_0' \geq -\frac{k}{u_0}[\alpha n_0 - E_0]$, fails, i.e.,*

$$\text{(5.39)} \qquad \exists \alpha \in \mathbb{R}^+ \quad s.t. \quad u_0' < -\frac{k}{u_0}[\alpha n_0 - E_0].$$

*Proof.* Instead of using the implicit formula of the flow map $r$ in (5.20) we introduce a dimensionless parameter $\tau$ such that $r$ may be rewritten in terms of the parameter $\tau \in \mathbb{R}^+$,

$$r = \frac{R}{2}[1 + \cosh(\tau + \tau_0(\alpha))],$$

$$t = \frac{R}{2Q}[\tau + \sinh(\tau + \tau_0(\alpha)) - \sinh(\tau_0(\alpha))],$$

where $Q$ and $R$ are given in (5.21), and $\tau_0$ is determined uniquely by

$$\cosh(\tau_0(\alpha)) = \frac{2\alpha}{R} - 1 = \frac{u_0^2}{ke_0}.$$

Thus the 'indicator' function $\Gamma$ is determined by

$$\Gamma(\alpha, t) = r_\alpha - r_\tau/t_\tau t_\alpha,$$

with

$$r_\alpha = \frac{R_\alpha}{2}[1 + \cosh(\tau + \tau_0(\alpha))] + \frac{R}{2}\sinh(\tau + \tau_0(\alpha))t_{0\alpha},$$

$$r_\tau = \frac{R}{2}\sinh(\tau + \tau_0(\alpha)),$$

$$t_\alpha = \frac{R}{2Q}[1 + \cosh(\tau + \tau_0(\alpha))],$$

$$t_\tau = \left(\frac{R}{2Q}\right)_\alpha [\tau + \sinh(\tau + \tau_0(\alpha)) - \sinh(\tau_0(\alpha))] + \frac{R}{2Q}[\cosh(\tau + \tau_0(\alpha)) - \cosh(\tau_0(\alpha))]\tau_{0\alpha}.$$

Expressed in terms of

$$A(\tau) := \frac{\sinh(\tau + \tau_0)}{1 + \cosh(\tau + \tau_0)},$$

we have

$$\Gamma(\alpha, t) = \frac{1 + \cosh(\tau + \tau_0)}{2}\Big\{R_\alpha + RA(\tau)\tau_{0\alpha}$$
$$- QA(\tau)\Big[(R/Q)_\alpha \left(\frac{\tau - \sinh(\tau_0)}{1 + \cosh(\tau + \tau_0)} + A(\tau)\right) + \frac{R}{Q}\left(1 - \frac{1 + \cosh(\tau_0)}{1 + \cosh(\tau + \tau_0)}\right)\tau_{0\alpha}\Big]\Big\}.$$

To show the breakdown of solutions it suffices to show that as $\tau$ becomes large $\Gamma$ becomes negative since $\Gamma(\alpha, 0) = 1 > 0$. Note that as $\tau \to \infty$ one has $A(\tau) \to 1$, and therefore the limit of the sum in the bracket $\{\cdots\}$ combined with $e_0 = E_0\alpha^2$ and $e_0' = \rho_0 = n_0\alpha^2$ becomes

$$\frac{R}{Q}Q_\alpha = \frac{R}{Q^2}\left[u_0 u_0' + k(\alpha n_0 - E_0)\right].$$

Hence the 'indicator' function $\Gamma$ would become negative for large time whenever the condition (5.39) is satisfied. Therefore there must be a finite $\tau_c$, also a finite time $t_c$, such that $\Gamma(\alpha, t_c)$ vanishes. At this time, $t_c$, the solution breaks down. $\square$

*Remark* 5.4. The above lower threshold for breakdown, $\nu = 2$, can also be verified in an alternative way as a particular case of the more general situation $\nu \geq 2$. Indeed, the new feature for $\nu \geq 2$, is that the velocity tends to a constant for large time t's, which is evident from the energy integral (5.30),

$$u^2 = (r')^2 = u_0^2 + \frac{2ke_0}{\nu - 1}(\alpha^{1-\nu} - r^{1-\nu}),$$

while noting that, say by (5.33), $r^{1-\nu} \to 0$ as $t \to \infty$. Hence, the velocity approaches the constant value

$$u(t) \sim Q(\alpha), \quad Q(\alpha) = \sqrt{u_0^2 + \frac{2ke_0}{\nu - 1}\alpha^{1-\nu}},$$

which in turn implies that

$$r \sim Q(\alpha)t.$$

Consequently, if the following critical condition fails,

$$u_0' \geq -\frac{k}{u_0}\left[\frac{n_0\alpha}{\nu - 1} - E_0\right],$$

i.e., if there exists an $\alpha \in \mathbb{R}^+$ such that $Q_\alpha < 0$, or explicitly, that

(5.40) $$\exists \alpha \in \mathbb{R}^+ \quad s.t. \quad u_0' < -\frac{k}{u_0}\left[\frac{n_0\alpha}{\nu - 1} - E_0\right],$$

then two particle paths must collide at large time. Observe that this critical condition for $\nu = 2$ coincides with (5.39).

We conclude with the general $\nu \geq 2$ case, discussing the upper threshold for the existence of global smooth solution.

**Theorem 5.9** (Global existence of smooth solutions for general $\nu \geq 2$ cases). *A global smooth solution of the Euler-Poisson equations (5.1)-(5.3) with $\nu \geq 2$ exists provided the initial data $(u_0, n_0)$ with $E_0 = \alpha^{-\nu}\int_0^\alpha n_0(\xi)\xi^\nu d\xi$ is prescribed such that for all $\alpha \in \mathbb{R}^+$*

(5.41) $$u_0' > -\frac{k}{u_0}\left[\frac{n_0\alpha^\nu h_\nu}{\nu - 1} - E_0\right].$$

*Here, $h_\nu$ is determined by*

(5.42) $$\frac{ku_0n_0\alpha}{(\nu - 1)^2}\int_0^{h_\nu}\frac{h_\nu(\alpha) - \eta}{[u_0^2 + \frac{2ke_0}{\nu-1}\eta]^{3/2}}(\alpha^{1-\nu} - \eta)^{\frac{\nu}{1-\nu}}d\eta \equiv 1.$$

*Remark* 5.5. Since, as we shall see below, $h_\nu < \alpha^{1-\nu}$, we conclude that the lower threshold (5.40) is indeed smaller than the upper threshold in (5.41).

*Proof.* Recalling the energy identity (5.30),

$$(r')^2 = u^2 = u_0^2 + \frac{2ke_0}{\nu - 1}(\alpha^{1-\nu} - r^{1-\nu}),$$

we have for $u_0 > 0$

$$\frac{dr}{\sqrt{u_0^2 + \frac{2ke_0}{\nu-1}(\alpha^{1-\nu} - r^{1-\nu})}} = dt.$$

Integration yields
$$\int_\alpha^{r(\alpha,t)} \frac{d\xi}{\sqrt{u_0^2 + \frac{2ke_0}{\nu-1}(\alpha^{1-\nu} - \xi^{1-\nu})}} = t.$$

Differentiating the above equality with respect to $\alpha$ leads to
$$\frac{\Gamma(\alpha,t)}{\sqrt{u_0^2 + \frac{2ke_0}{\nu-1}(\alpha^{1-\nu} - r^{1-\nu})}} = B(r,t)$$

where
$$B(\alpha,t) := \int_\alpha^r \frac{u_0 u_0' + \frac{k}{\nu-1}\rho_0 \alpha^{1-\nu} - ke_0 \alpha^{-\nu} - \frac{k}{\nu-1}\rho_0 \xi^{1-\nu}}{[u_0^2 + \frac{2ke_0}{\nu-1}(\alpha^{1-\nu} - \xi^{1-\nu})]^{3/2}} d\xi + \frac{1}{u_0},$$

which can be rewritten in terms of $M_\nu := (\nu-1)\frac{kE_0 - u_0 u_0'}{kn_0 \alpha^\nu}$ and $\eta := \alpha^{1-\nu} - \xi^{1-\nu}$ as
$$B(\alpha,t) = \frac{kn_0 \alpha}{(\nu-1)^2} \int_0^{\alpha^{1-\nu} - r^{1-\nu}} \frac{\eta - M_\nu}{[u_0^2 + \frac{2ke_0}{\nu-1}\eta]^{3/2}} (\alpha^{1-\nu} - \eta)^{\frac{\nu}{1-\nu}} d\eta + \frac{1}{u_0}.$$

In remains to show the positivity of $B$ for all $t > 0$ provided (5.41) holds. If $M_\nu \leq 0$, then it is easy to see $B > 0$ for all $t > 0$. We now consider the case $M_\nu > 0$ by checking the positivity of the possible minimum of $B$. Note that $e_{0\alpha} = \rho_0 = n_0 \alpha^\nu$ and $e_0 = E_0 \alpha^\nu$. A straightforward calculation involving (5.30) gives
$$\frac{dB}{dt} = \frac{kn_0 \alpha}{\nu-1} \frac{\alpha^{1-\nu} - r^{1-\nu} - M_\nu}{u_0^2 + \frac{2ke_0}{\nu-1}(\alpha^{1-\nu} - r^{1-\nu})}.$$

¿From the monotonicity of the flow map $\frac{dr}{dt} > 0$, it follows that there exists a time $t^*$ such that
$$\frac{dB}{dt}|_{t=t^*} = 0, \quad \frac{dB}{dt}(t - t^*) > 0 \quad \text{for} \quad t \neq t^*,$$

and at this time $r(\alpha, t^*)^{1-\nu} = \alpha^{1-\nu} - M_\nu$. Therefore we have
$$B(\alpha,t) \geq B(\alpha, t^*) = \frac{kn_0 \alpha}{(\nu-1)^2} \int_0^{M_\nu} \frac{\eta - M_\nu}{[u_0^2 + \frac{2ke_0}{\nu-1}\eta]^{3/2}} (\alpha^{1-\nu} - \eta)^{\frac{\nu}{1-\nu}} d\eta + \frac{1}{u_0}$$
$$> \frac{kn_0 \alpha}{(\nu-1)^2} \int_0^{h_\nu} \frac{\eta - h_\nu(\alpha)}{[u_0^2 + \frac{2ke_0}{\nu-1}\eta]^{3/2}} (\alpha^{1-\nu} - \eta)^{\frac{\nu}{1-\nu}} d\eta + \frac{1}{u_0} = 0,$$

provided for all $\alpha \in \mathbb{R}^+$, $M_\nu(\alpha) < h_\nu(\alpha) < \alpha^{1-\nu}$, i.e., (5.41), with $h^\nu$ defined in (5.42). Hence $\Gamma(\alpha, t)$ remains positive for all $t > 0$ once the initial data remain above the upper threshold (5.41). □

The above upper and lower thresholds for the cases $\nu \geq 1$ confirm the existence of an intermediate critical threshold, through we do not provide the explicit form of the critical threshold. One case in which we can precisely compute the critical threshold is the 4-dimensional ($\nu = 3$) isotropic case, which is given in the following

**Theorem 5.10** (Critical threshold for 4-dimensional model ($\nu = 3$)). *The global smooth solution of Euler-Poisson equations (5.1)-(5.3) with $\nu = 3$ exists if and only if*

(5.43) $$(\alpha u_0' - u_0)^2 < 4\alpha(\frac{n_0 \alpha}{2} - E_0), \quad \forall \alpha > 0 \quad \text{and} \quad u_0'(0) \geq 0.$$

*In this case the velocity is given by*

$$u(r,t) = \frac{\alpha u_0 + [u_0^2 + k\alpha E_0]t}{\sqrt{\alpha^2 + 2\alpha u_0 t + [u_0^2 + k\alpha E_0]t^2}} \to \sqrt{u_0^2 + k\alpha E_0} \quad \text{as} \quad t \uparrow \infty, \tag{5.44}$$

*and the density is given by*

(5.45)
$$n(r(\alpha,t),t) = \frac{n_0(\alpha)\alpha^3}{[\alpha^2 + 2\alpha u_0 t + (u_0^2 + kE_0\alpha)t^2][\alpha + (u_0 + \alpha u_0')t + (u_0 u_0' - kE_0 + \tfrac{k}{2}n_0\alpha)t^2]}.$$

*Proof.* As argued before it is sufficient and necessary to show that the threshold condition (5.43) ensures the positivity of the indicator function $\Gamma$ for all $t > 0$. Let us first solve the flow map equation (5.11) with $\nu = 3$, i.e.,

$$r'' = ke_0 r^{-3}, \quad r(0) = \alpha, \quad r'(0) = u_0.$$

Its energy integral is

$$[r']^2 = u_0^2 + ke_0\alpha^{-2} - ke_0 r^{-2}$$
$$= u_0^2 + ke_0\alpha^{-2} - rr'',$$

where $r'' = ke_0 r^{-3}$ has been used in the last equalities. Rewriting this relation leads to

$$\frac{1}{2}[r^2]'' = [r']^2 + rr'' = u_0^2 + ke_0\alpha^{-2},$$

and integration twice gives

$$r^2 = \alpha^2 + 2\alpha u_0 t + [u_0^2 + ke_0\alpha^{-2}]t^2.$$

Hence

$$r = \sqrt{\alpha^2 + 2\alpha u_0 t + [u_0^2 + ke_0\alpha^{-2}]t^2}, \tag{5.46}$$

which as $u = \frac{dr}{dt}$ and $e_0 = E_0\alpha^3$ gives the velocity (5.44). The corresponding 'indicator' function is

$$\Gamma(\alpha,t) = \frac{\partial r}{\partial \alpha} = \frac{\alpha + [\alpha u_0]'t + \tfrac{1}{2}[u_0^2 + k\alpha E_0]'t^2}{\sqrt{\alpha^2 + 2\alpha u_0 t + [u_0^2 + k\alpha E_0]t^2}}$$
$$= \frac{\alpha + [u_0 + \alpha u_0']t + [u_0 u_0' - kE_0 + \tfrac{kn_0\alpha}{2}]t^2}{\sqrt{\alpha^2 + 2\alpha u_0 t + [u_0^2 + k\alpha E_0]t^2}}.$$

Note that at the origin $\alpha = 0$, $\Gamma(0,t) = 1 + u_0'(0)t$. These explicit formula imply that $\Gamma(\alpha,t)$ remains positive for all $t > 0$ once

$$[u_0 + \alpha u_0']^2 < 4\alpha(u_0 u_0' - kE_0 + \frac{kn_0\alpha}{2}), \quad \forall \alpha > 0 \quad \text{and} \quad u_0'(0) \geq 0.$$

This is equivalent to (5.43). The solution (5.45) follows from the above explicit expression of $\Gamma(\alpha,t)$ when recalling the general formula (5.13). $\square$

The above critical threshold result enables us to claim the following

**Corollary 5.11** (Breakdown of smooth solutions for the case $\nu = 3$). *A solution of the Euler-Poisson equations (5.1)-(5.3) with $\nu = 3$ blows up in finite time if and only if condition (5.43) fails, i.e.,*

$$\exists \alpha \in \mathbb{R}^+ \text{s.t.} \quad (\alpha u_0' - u_0)^2 \geq 4k\alpha(\frac{n_0\alpha}{2} - E_0).$$

*In this case, $n(r,t)$ and $u_r(r,t)$ become infinite as $t \uparrow T$, where the blow-up time is given explicitly by*

$$t_c := 2/\sup\left\{-u_0' - \frac{u_0}{\alpha} + \frac{1}{\alpha}\sqrt{(\alpha u_0' - u_0)^2 - 4k\alpha(\frac{n_0\alpha}{2} - E_0)}\right\}.$$

## 6. Appendix

*Proof of the maximum principle:*

To show that $u_t + f(x,t)u_x = a(x,t)u_{xx}$ satisfies a maximum principle we find a function that satisfies:

$$F''(x) = F(x)/(2b(x))$$

where $b(x) > a(x,t)$, and $b(x) \geq cx^2 + 1$. We then consider the function $z(x,t) = e^t(F(x) + d)$ and following [13] we prove the maximum principle.

We make use of general results about the solutions of $F''(x) = F(x)/(2b(x))$ in the proof. Because $b(x) \geq cx^2 + 1, c > 0$, it is easy to show [4, p. 106, question 35] that the solution of the ODE must look like a solution of $F''(x) = 0$ for large enough $x$; i.e. $F(x)$ looks like $r_{\pm}x + s_{\pm}$. It is easy to see that if we let $F(x_0) = 1$ and $F'(x_0) = 0$, then the solution of the ODE will be concave up near $x_0$. As $F''(x)$ is positive at $x_0$, we find that to the left of $x_0$, the function $F(x)$ must be positive and its derivative must be negative. Similarly, to the right of $x_0$ $F(x)$ must be positive and its derivative must be positive. Thus $F(x) \geq 1$ for all $x$. Thus, $F''(x)$ is also positive everywhere. That implies that $F'(x)$ increase from some negative value to some positive one. Since asymptotically $F(x)$ is linear, we find that $F'(x)$ must tend to some constant value as $x$ tends to $\pm\infty$. In what follows, we shall point out $F(x)$'s dependence on $x_0$ by using the notation $F_{x_0}$ to denote the solution of our ODE with initial data given at $x_0$.

We note that if we define $G_{x_0,d} = F_{x_0}(x) + d$, then for all positive $d$, $G_{x_0,d}$ satisfies the inequality $G''_{x_0,d}(x) \leq G_{x_0,d}/(2b(x))$. As $F'_{x_0}$ is bounded, as $d$ increase $|G'_{x_0,d}/G_{x_0,d}|$ tends to zero uniformly in $d$.

Consider the function $z(x,t) = G_{x_0,d}(x)e^t$. We find that $z_t = z$ and $z_{xx} = G''_{x_0,d}(x)e^t$. Thus,

$$\begin{aligned}z_t + f(x,t)z_x &= G_{x_0,d}(x)e^t + f(x,t)G'_{x_0,d}(x)e^t \\ &= e^t\left(G_{x_0,d}(x)/2 + b(x)G''_{x_0,d}(x) + f(x,t)G'_{x_0,d}(x)\right) \\ &\geq a(x,t)z_{xx} + e^t\left(G_{x_0,d}(x)/2 + f(x,t)G'_{x_0,d}(x)\right).\end{aligned}$$

If we make $d$ large enough we can make $G$ as much larger than $G'$ as we please. Thus, we find that for sufficiently large $d$, $z_t + z_x f(x,t) > a(x,t)z_{xx}$. Moreover, $z(x,0)$ tends to infinity linearly in $x$ and exponentially in $t$.

We note that we can use $e^{ct}$ rather than $e^t$ by letting $F$ solve the equation $F''(x) = cF(x)/(2b(x))$. This is what allows us to state that $E(t)$ may be exponential and need not be sub-exponential.

Now we modify one of the standard proofs of the maximum principle for the heat equation [13, pp. 216-218]. We consider $w(x,t) = u(x,t) - \epsilon z(x,t)$ where $\epsilon > 0$. Clearly, $w$ satisfies $w_t + f(x,t)w_x < a(x,t)w_{xx}$. If we consider this equation on a finite interval, $[x_1, x_2]$, then $w$ satisfies the maximum principle:

$$w(x,t) \leq \max(w(x_1,t), w(x_2,t), \max_{x \in (x_1,x_2)} w(x,0)).$$

As $z(x,t) > 0$, we find that $\max_{x \in (x_1, x_2)} w(x,0) < \sup_x u(x,0)$. Also, as $z(x,t)$ increase linearly in space and exponentially in time (with any desired exponent) and, by assumption, $u(x,t)$ increase more slowly, it is clear that for any $x_0$, for any $x$ of sufficiently large magnitude, $w(x,t) < \sup_x u(x,0)$. In fact it is not necessary that $u(x,t)$ be strictly sub-linear. It is sufficient that it be sub-linear on an infinite sequence of points for which $\pm\infty$ are limit points. If $|u(x,t)| < D(x)e^{\alpha t}$ for all $t$ and an infinite sequence of values, $\{x_i\}$, that has as limit points $\pm\infty$, then we will be able to find an infinite sequence of points, $\tilde{x}_i$, for which $w(x,t) < \sup_x u(x,0)$.

Putting all of this together, we find that for any $\epsilon$ if the magnitude of $x_1$ and $x_2$ are sufficiently large and $x_1$ and $x_2$ belong to the sequence of points on which $u(x,t)$ is sublinear, then $w(x_1,t), w(x_2,t) < \sup_x u(x,0)$. Finally we note that as $\epsilon \to 0$, $w(x_0, t) \to u(x_0, t)$. Thus, we find that $u(x_0, t) \leq \sup u(x, 0)$. As our argument in no way depends on $x_0$, we find that $\sup_x u(x,t) \leq \sup_x u(x,0)$. A simple corollary of this is that if all the hypotheses above are met, then if $u_t + f(x,t)u_x = c + a(x,t)u_{xx}$, then $s(x,t) = u(x,t) - ct$ satisfies the maximum principle. Thus, $u(x,t) \leq \sup_x u(x,0) + ct$. Finally, by also considering $-u$ we find that $\inf_x u(x,0) + ct \leq u(x,t) \leq \sup_x u(x,0) + ct$.

*Remark* 6.1. We see that if it is known that $|u(x,t)| \leq D(x)e^{\alpha t}$ for all $t$ and for an infinite sequence of $x$'s that run to $\pm\infty$, then $u(x,t)$ satisfies a maximum principle. If one knows that the growth of $u(x,t)$ in time is only exponential, then in order for $u$ to fail to satisfy a maximum principle, it is necessary that $u(x,t)$ have growth in $x$ that is faster than any sub-linear function. In particular, $u(x,t) > Mx/ln(x)$ for all sufficiently large $x$.

*Remark* 6.2. We show that having some condition on $f(x,t)$ is necessary. Consider the equation $u_t = b(x)u_{xx}$. If we let $v = u_x$, then we find that $v_t = (b(x)v_x)_x$ We have already shown that if $b(x) > x^{2+\epsilon} + 1$, then the equation $F''b(x) = F$ has a solution that grows linearly at $\pm\infty$ and whose derivative, $F'(x)$ increase from some value at $-\infty$ to some value at $\infty$. Let $H(x) = F'(x)$. We find that $(b(x)H'(x)) = H(x)$. Thus, $e^t H(x)$ is a solution of the equation for $v$. This is a bounded solution of the equation that satisfies neither a maximum nor a minimum principle. Of course, we can rewrite the PDE for $v$ in the form $v_t - b'(x)v_x = b(x)v_{xx}$. Thus we see that it is imperative that some conditions be placed on $f(x,t)$. Clearly this "hyperbolic term" *can* destabilize the parabolic PDE.

If one specializes to $b(x)$ which are even, then it is easy to say more. For such $b(x)$ it is easy to see that the solution of $b(x)F'' = F$ with initial data $F(0) = 0, F'(0) = a, a \neq 0$ is odd, linear at infinity, and $F'(x)$ is even and tends to a nonzero constant, $k$, at $\pm\infty$. $F'(x)$ will always be greater than or equal to $a$. If $b(x) = \cosh(x)$, then (using more results on the asymptotic behavior of ODEs) we also find that $F(x) \to kx + l$ exponentially fast at $\pm\infty$. Consider the function $e^t F(x) - kx - l$. It is initially bounded (it even tends to 0 at $\pm\infty$), it is a solution of $u(x,t) = b(x)u_{xx}$, and it does not satisfy a maximum principle. We note that a function that solves the heat equation cannot blow up in this fashion. The bounds on the solution and the fact that the solution was initially bounded would be enough to guarantee that the solution remained bounded.

## ACKNOWLEDGMENTS

Research was supported in part by ONR Grant No. N00014-91-J-1076 (ET) and by NSF grant #DMS97-06827 (ET, HL). Additional support was provided by a Jerusalem College of Technology Presidential Research Grant (SE). H. Liu wants to thank Heinz-Otto Kreiss for his many generous discussions.

## References

[1] U. M. Ascher, Peter A. Markowich, P. Pietra, and C. Schmeiser, *A phase plane analysis of transonic solutions for the hydrodynamic semiconductor model*, Math. Models Methods Appl. Sci. **1** (1991), 347–376.

[2] U. Brauer, A. Rendal, and O. Reula, *The cosmic no-hair theorem and the non-linear stability of homogeneous Newtonian cosmological models*, Class. Quantum Grav. **11** (1994), 2283-2296.

[3] C. Cercignani, *The Boltzmann Equations and Its Applications*, Springer- Verlag, New York, 1988.

[4] E.A. Coddington and N. Levinson, *Theory of Ordinary Differential Equations*, Robert E. Krieger Publishing Company, Malabar, Florida, 1987.

[5] G.-Q. Chen and D. Wang, *Convergence of shock capturing scheme for the compressible Euler-Poisson equations*, Comm. Math. Phys. **179** (1996), 333-364.

[6] J. Dolbeault and G. Rein, *Time-dependent rescalings and Lyapunov functionals for the Vlasov-Poisson and Euler- Poisson systems, and for related models of kinetic equations, fluid dynamics and quantum physics*, Mathematical Models and methods in Applied Sciences.

[7] S. Engelberg, *Formation of singularities in the Euler-Poisson equations*, Physica D **98** (1996), 67-74.

[8] P. Gamblin, *Solution rgulire temps petit pour l'quation d'Euler-Poisson*, Comm. Partial Differential Equations **18** (1993), 731–745.

[9] I. Gasser, C-K Lin and P.A. Markowich, *A review of dispersive limits of (non)linear Schrödinger-type equations*, preprint, 2000.

[10] Y. Guo, *Smooth irrotational flows in the large to the Euler-Poisson system in $\mathbb{R}^{3+1}$*, Comm. math. Phys. **195** (1998), 249-265.

[11] D. Holm, S. F. Johnson, and K.E. Lonngren, *Expansion of a cold ion cloud*, Appl. Phys. Lett. **38** (1981), 519-521.

[12] J.D. Jackson, *Classical Electrodynamics*, 2nd ed., Wiley, New York, 1975.

[13] F. John, *Partial Differential Equations, Fourth Edition*, Springer-Verlag, New York 1982.

[14] S. Junca and M. Rascle, *Relaxation of the isothermal Euler-Poisson system to the drift-diffusion equations*, Quart. Appl. Math. **58** (2000), 511-521.

[15] T. Luo, R. Natalini, and Z. Xin, *Large time behavior of the solutions to a hydrodynamic model for semiconductors*, SIAM J. Appl. Math. **59** (1999), 810–830.

[16] T. Makino, *On a local existence theorem for the evolution of gaseous stars,* Patterns and Waves (eds. T. Nishida, M. Mimura and H. Fujii), North-Holland/Kinokuniya, 1986, 459-479.

[17] P.A. Markowich, *A non-isentropic Euler-Poisson Model for a Collisionless Plasma*, Math. Methods Appl. Sci. **16** (1993), 409-442.

[18] T. Makino, and B. Perthame, *Sur les solutions a symetrie spherique de l'equation d'Euler-Poisson pour l'etoiles gazeuses*, Japan J. Appl. math., **7** (1990), 165-170.

[19] P. Marcati and R. Natalini, *Weak solutions to a hydrodynamic model for semiconductors and relaxation to the drift- diffusion equation*, Arch. Rat. Mech. Anal. **129** (1995), 129- 145.

[20] P.A. Markowich, C. Ringhofer, and C. Schmeiser, *Semiconductor Equations*, Springer, Berlin, Heidelberg, New York, 1990.

[21] T. Makino and S. Ukai , *Sur l'existence des solutions locales de l'quation d'Euler-Poisson pour l'volution d'toiles gazeuses*, J. Math. Kyoto Univ. **27** (1987), 387–399.

[22] B. Perthame, *Nonexistence of global solutions to the Euler-Poisson equations for repulsive forces*, Japan J. Appl. Math. **7** (1990), 363-367.

[23] P. Rosenau, *Evolution and breaking of ion-acoustic waves*, Phys. Fluids **31** (6) (1988), 1317- 1319.

[24] D. Wang, *Global solutions and relaxation limits of Euler-Poisson equations*, to appear in Z. Angew. Math. Phys.

[25] D. Wang, *Global solutions to the equations of viscous gas flows*, to appear in Pro. Royal Soc. Edinburgh: Section A.

[26] D. Wang *Global solutions to the Euler-Poisson equations of two-carrier types in one-dimension*, Z. Angew. Math. Phys. **48** (1997), 680-693.

[27] D. Wang, G.-Q. Chen, *Formation of singularities in compressible Euler-Poisson fluids with heat diffusion and damping relaxation*, J. Diff. Eqs. **144** (1998), 44- 65.

Electronics Department, Jerusalem College of Technology—Machon Lev, P.O.B. 16031, Jerusalem, Israel
 *E-mail address*: `shlomoe@optics.jct.ac.il`

UCLA, Department of Mathematics, Los Angeles, CA 90095-1555, USA
 *E-mail address*: `hliu@math.ucla.edu`

UCLA, Department of Mathematics, Los Angeles, CA 90095-1555, USA
 *E-mail address*: `tadmor@math.ucla.edu`